\documentclass[dvips]{amsart}
\usepackage[dvips]{graphics}
\usepackage{amscd,amssymb,amsxtra,latexsym,psfrag,epsfig}
\newcommand{\vanish}[1]{}
\newcommand{\hz}{\hat{0}}
\newcommand{\ho}{\hat{1}}
\newcommand{\Ff}{{\mathcal F}}
\newcommand{\tensor}{\otimes}

\newcommand{\D}{\Delta}
\newcommand{\Dhat}{\hat \Delta}
\newcommand{\Fhat}{\hat{\mathcal F}}
\newcommand{\into}{\hookrightarrow}

\newcommand{\ab}{\av\bv}
\newcommand{\av}{{\bf a}}
\newcommand{\bv}{{\bf b}}
\newcommand{\cd}{\cv\dv}

\newcommand{\cv}{{\bf c}}
\newcommand{\dv}{{\bf d}}
\newcommand{\ev}{e}

\newcommand{\id}{\operatorname{id}}
\newcommand{\wt}{\operatorname{wt}}
\newenvironment{proof_}[1]{\noindent {\bf #1}}{{\qed}}
\theoremstyle{plain}
\newtheorem{theorem}{Theorem}
\newtheorem*{theorem*}{Theorem}
\newtheorem{proposition}{Proposition}
\newtheorem{conjecture}{Conjecture}
\newtheorem{corollary}{Corollary}

\newtheorem{lemma}{Lemma}[section]

\theoremstyle{remark}
\newtheorem*{remark}{Remark}
\newtheorem{example}{Example}
\newtheorem{exercise}{Exercise}

\newcommand{\marcelo}{aguiar00:_infin_hopf}
\newcommand{\marcelocd}{aguiar02:_infin_hopf}
\newcommand{\bk}{bayer91}
\newcommand{\bb}{bayer85:_gener_dehn_sommer_euler}
\newcommand{\bn}{billera00:_noncom}
\newcommand{\be}{billera00:_monot}
\newcommand{\ber}{billera97}
\newcommand{\ehren}{ehrenborg02:_inequal}
\newcommand{\ehrenborg}{ehrenborg01:_euler}
\newcommand{\emahajan}{ehrenborg98:_maxim}
\newcommand{\ereaddyr}{ehrenborg96}
\newcommand{\efox}{ehrenborg:_inequal}
\newcommand{\ereaddy}{ehrenborg98:_coprod}
\newcommand{\joni}{joni79:_coalg}
\newcommand{\kalai}{kalai88}
\newcommand{\purtill}{purtill93:_andre}
\newcommand{\reading}{reading02:_bases_euler}
\newcommand{\stanleycd}{stanley94:_flag}
\newcommand{\stanleylog}{stanley89:_log}
\newcommand{\stenson}{stenson02:_relat}
\newcommand{\sundaram}{sundaram94:_cohen_macaul}
\newcommand{\sund}{sundaram95}

\begin{document}
\title{The {\bf cd}-index of the Boolean lattice}
\author{Swapneel Mahajan}

\address{Department of Mathematics\\
Cornell University\\
Ithaca, NY 14853}
\email{swapneel@math.cornell.edu}
%\date{\today}
%\date{September 20, 1999}


\begin{abstract}

We study some properties 
of the {\bf cd}-index of the Boolean lattice.
They are extremely similar 
to the properties of the {\ab}-index,
or equivalently, the flag $h$-vector
of the Boolean lattice and hence may
be viewed as their {\bf cd}-analogues. 
We define a different algebra structure on the 
polynomial algebra $k \langle \cv, \dv \rangle$
and give a 
derivation on this algebra.
It is of significance for the Boolean lattice
and forms our main tool.
Using similar methods, we also prove some results 
for the {\bf cd}-index of the cubical lattice.
We show that the Dehn-Sommerville relations for
the flag $f$-vector of an Eulerian poset are equivalent
to certain simple identities that exist in our 
algebra.

\end{abstract}

\maketitle

\section{Introduction}

The {\bf cd}-index is a non-commutative polynomial 
in the variables {\bf c} and {\bf d} which efficiently 
encodes the flag
$f$-vector (equivalently the flag $h$-vector) of an
Eulerian poset. 
The flag $h$-vector and the {\bf cd}-index are mysterious
objects with many interesting properties. 
It is true, for example, that
the {\bf cd}-index (and hence the {\bf ab}-index)
of the face lattice of a convex polytope 
is a polynomial with positive integer coefficients.
We refer the
reader to the basic paper of Stanley~\cite{\stanleycd}.
For more recent references, 
see~\cite{\marcelocd,\be,\ber,\bn,\ehrenborg,\efox,\ereaddy}.
In this section, we first review the basic definitions
and then motivate the problem that we study.

\subsection{Background}

Let $P$ be a graded partially ordered set (poset) of rank $n+1$.  
For $S$, a subset of $[n]$,
%$\{1,\ldots,n\}$, 
let $f_{S}$ be the number of chains (flags) in $P$
that have elements on the ranks in $S$. 
These $2^n$ numbers
constitute the flag $f$-vector of the poset $P$. 
The flag $h$-vector
is defined by the relation
$h_{S} = \sum_{T \subseteq S} (-1)^{|S-T|} f_{T}$.
The standard way to encode the flag $h$-vector 
is to express it as a non-commutative polynomial
in the variables $\av$ and $\bv$ as follows.
Define a monomial
$u_S = u_1 u_2 \ldots u_n$ by
$$u_i = \left\{\begin{array}{c c}
             \av, & i \not \in S \\
             \bv, & i \in S,
             \end{array}\right.$$
\noindent
and let 
$\Psi_P(\av,\bv) = \sum_{S \subseteq [n]} h_S u_S$.
The polynomial $\Psi_P(\av,\bv)$ is called 
the {\bf ab}-index of $P$.

A poset $P$ is called Eulerian 
if for all $x \leq y$ in $P$,
we have 
$\mu(x,y) = (-1)^{\rho(x,y)}$,
where $\mu$ denotes the Mobius function of 
the interval $(x,y)$ of $P$
and where $\rho(x,y) = \rho(y) - \rho(x)$.
Here $\rho$ is the rank function of $P$.
An important example of an Eulerian poset is the face
lattice of a convex polytope.

\begin{theorem*}
If $P$ is Eulerian then 
$\Psi_P(\av,\bv)$ can be written uniquely as a polynomial
$\Psi_P(\cv,\dv)$ in the non-commuting variables 
$\cv = \av+\bv$ and $\dv = \av\bv+\bv\av$.
\end{theorem*}

This fact was noticed by Fine and
proved by Bayer and Klapper,
see~\cite{\bk} or~\cite[Theorem 1.1]{\stanleycd}.
It is equivalent to the statement that
the linear relations satisfied by the flag $h$-vector
of an Eulerian poset are precisely the
generalised Dehn-Sommerville equations, also known as
the Bayer-Billera relations~\cite{\bb}.

The polynomial $\Psi_P(\cv,\dv)$ 
is called the {\bf cd}-index of $P$.
Henceforth we will just refer to 
it as $\Psi(P)$.
Note that the degree of $\Psi(P)$ 
is one lower than the rank of the poset $P$. 
(The variables $\cv$ and $\dv$ are assigned degrees $1$ and
$2$ respectively.)
It also follows from the definition that
for a poset $P$ of rank $n+1$,
the coefficient of $\cv^n$ in $\Psi(P)$ is always $1$.

\vanish{
It gives an explicit basis for the
generalised Dehn-Sommerville equations, also known as
the Bayer-Billera relations. 
For a poset $P$, we denote its
{\bf cd}-index by
$\Psi(P)$. 
It is a homogeneous polynomial in the non-commuting variables 
{\bf c} and {\bf d}, 
}

\begin{figure} 
$$  \begin{array}{r c l}
\Psi(B_0) & = & \ev \\
\Psi(B_1) & = & 1 \\
\Psi(B_2) & = & \cv \\
\Psi(B_3) & = & \cv^2 + \dv \\
\Psi(B_4) & = & \cv^3 + 2 \cv \dv + 2 \dv \cv \\
\Psi(B_5) & = & \cv^4 + 3 (\cv^2 \dv + \dv \cv^2) +
                5 \cv \dv \cv + 4 \dv^2
     \end{array} $$
\caption{The {\bf cd}-index of the Boolean lattice for small ranks.}
\label{f:bl}
\end{figure}

\subsection{The Boolean lattice}

The Boolean lattice of rank $n+1$, 
which we denote by $B_{n+1}$, 
is the poset of all subsets of the set $[n+1]$
%$\{1,2,\ldots,n+1\}$
ordered by inclusion.
It is same as the face lattice of the $n$-dimensional simplex.
Hence it is an Eulerian poset and has a 
{\bf cd}-index.
Since the simplex is the simplest polytope,
the Boolean lattice has a special role to play
in the class of Eulerian posets.
For example, among all Eulerian posets of rank $n+1$,
the Boolean lattice $B_{n+1}$ has 
the smallest $\cd$-index coefficient-wise.
This was a conjecture of Stanley 
which was proved by Billera and Ehrenborg~\cite{\be}. 

In this paper,
we will be studying the Boolean lattice,
not in relation to other posets,
but rather as an object by itself.
The flag $h$-vector (or in other words, the {\bf ab}-index)
of the Boolean lattice 
displays many remarkable patterns.
These have been studied in detail in~\cite{\emahajan}.
The goal of this paper is to establish
the {\bf cd}-analogues of these properties
for the Boolean lattice.

The intuitive reason why
properties of the $\ab$-index
carry over to the $\cd$-index
is explained in Appendix~\ref{s:connection}.
Unfortunately,
the methods used to study the two problems
are totally different.
The coefficients of the $\ab$-index of the Boolean lattice
are related to the descent statistic of permutations
and they can be studied effectively
using that interpretation~\cite{\emahajan}.
Though there are similar interpretations
for the coefficients of the $\cd$-index 
(see item (1) in Section~\ref{ss:questions}),
we do not know how to use them.
Instead, our method is based on an algebraic study of the
polynomial algebra 
${\mathcal F} = k \langle \cv,\dv \rangle$ 
in the non-commuting
variables {\bf c} and {\bf d}. 

The {\bf cd}-index of the Boolean lattice for small ranks
is shown in Figure~\ref{f:bl}.
The letter $\ev$ is a formal symbol that is added to $\Ff$
in degree $-1$, see Section~\ref{ss:ext}.

\subsection{Questions (and partial answers)}
\label{ss:questions}

For $v$, a {\bf cd}-monomial of degree $n$,
let $\beta(v)$ be the coefficient of 
$v$ in $\Psi(B_{n+1})$.
Our primary objective is to understand $\beta$.

\begin{enumerate}
\item[(1)] 
We know that $\Psi(B_{n})$ is a polynomial whose coefficients
are positive integers. 
What do these numbers count? 

The {\bf cd}-index of $B_n$ is a refined enumeration
of Andr\'{e} permutations~\cite{\purtill}. 
Similarly, it is also a refined enumeration
of simsun permutations,
first defined by Simion and Sundaram~\cite{\sundaram,\sund}.
These permutations seem ad~hoc and
it is not clear how to use them to 
answer some of the questions asked below. 
More recently,
there has been an interpretation 
involving the peak statistic of permutations,
which may prove more useful.

\item[(2)] What can be said about the equalities
satisfied by the values of $\beta$ ?

It is known that
$\beta(v) = \beta(v^*)$,
where $v^*$ is the monomial $v$
written in reverse order.
For small ranks,
this can be seen from the data in Figure~\ref{f:bl}.
Hence one may ask:
Are there {\bf cd}-monomials $u$ and $v$
such that $u \not = v^*$, but yet
$\beta(u) = \beta(v)$?
Lemma~\ref{l:identities}
provides a partial answer to this question.
Lemma~\ref{l:switch} 
gives a slightly more general answer
that involves a different product $\centerdot$ on
$k \langle \cv,\dv \rangle$.
Also see Corollary~\ref{c:identity} in Section~\ref{s:ln}.
It is stated using a different notation involving lists.

\item[(3)] Among all {\bf cd}-monomials of a given degree, 
which {\bf cd}-monomial has the largest $\beta$ value?

We treat this problem in Section~\ref{s:maxima}.
Theorem~\ref{t:maxima} gives a simple and complete answer.
Though the main idea of the proof is simple,
we have to rely on two facts that are
stated as exercises.
This makes the proof a little unsatisfactory.

\item[(4)] What can be said about inequalities in general?

A simple and striking inequality 
is provided by Lemma~\ref{l:cc_d}.
It says that replacing an occurrence 
of $\cv^2$ by $\dv$ in a $\cd$-monomial
increases its $\beta$ value.
This is the first step for solving the question 
that was raised in item (3).
Ehrenborg~\cite{\ehren} has shown recently that
this inequality holds for the $\cd$-index of any polytope.

There are two types of inequalities
that we study in detail.
The first type occur as portions of 
reverse unimodal sequences
(Section~\ref{s:unimodal})
and the second type
are the balance inequalities 
(Section~\ref{s:balance}).
The motivation for the latter comes from a conjecture
of Gessel about the {\bf ab}-index of the Boolean lattice,
which was proved in~\cite{\emahajan}. 
We propose a {\cd}-analogue 
to this conjecture; see Conjecture~\ref{c:balance}
in Section~\ref{s:balance}.
There is plenty of evidence as to why
this conjecture should be true.
Theorems~\ref{t:adjoining_balance} and~\ref{t:suff_cond} are 
important results in this direction.

\item[(5)] What can be said about formulas for the $\beta$ values ?

We show that
$\beta(\cv^i \dv \cv^j) = {\binom {i+j+2}{i+1}} -1$
and
$\beta(\dv^n) = \frac{1}{2^n} E_{2n+1}$,
where $E_{2n+1}$ are the Euler or tangent numbers
(Example~\ref{eg:compute} and Lemma~\ref{l:deuler} respectively).
Further data suggests that studying exact values in detail
may be very interesting.
For example,
many of the $\beta$ values are divisible by 1001;
see Section~\ref{s:conc_rmks}.
However, the thrust of this paper is on studying inequalities.

\end{enumerate}

\subsection{Organisation of the paper}

We begin with the study of the polynomial algebra 
${\mathcal F} = k \langle \cv,\dv \rangle$ in
Section~\ref{s:algebra}.
Following Ehrenborg-Readdy~\cite{\ereaddy},
we define a coproduct $\Delta$ and a derivation $G$
on the algebra ${\mathcal F} = k \langle \cv,\dv \rangle$.
We first modify $\Ff$ to $\Fhat$ 
by adding a piece of degree $-1$ and then extend
$\Delta$ and $G$ to
$\Dhat$ and $\hat G$ respectively.
The main result of Section~\ref{s:algebra} is that
the extended map ${\hat G}$ is a coderivation on $\Fhat$.
The connection with the Boolean lattice
is provided by the identity
${\hat G}(\Psi(B_n)) = \Psi(B_{n+1})$.
The map $\hat G$ is our main tool for answering the questions in
Section~\ref{ss:questions}.

In Section~\ref{s:duality}, we dualise the maps
$\Dhat$ and $\hat G$ to get respectively
a product
(denoted by $\centerdot$) on $\Fhat$ 
and a derivation $S$, with the property
$S(v) = S(\beta(v))$.
We write down explicit formulas 
for the product $\centerdot$ 
and the derivation $S$.
By way of motivation, 
we provide quick applications to the Boolean lattice.
The explicit product $\centerdot$ has also appeared independently
in the work of Stenson and Reading; 
see in particular~\cite[Theorem 11]{\stenson} and
~\cite[Proposition 21]{\reading}.
The preprints are available on their respective homepages.
I thank Ehrenborg for pointing these references to me.

In Section~\ref{s:ln}, 
we restate all earlier results in
an alternate notation 
for {\cd}-monomials that involves lists.
This notation is quite natural
and easy to work with.
In Sections~\ref{s:unimodal},\ \ref{s:maxima} and \ref{s:balance},
we use the tools developed in earlier sections
to target two specific problems, namely
those raised in items (3) and (4) in Section~\ref{ss:questions}.
Throughout these sections,
we work with the list notation.
Wherever convenient, we also state our results in terms of 
the original notation of monomials.

There are four appendices.
In Appendix~\ref{s:connection},
we give some connection between the 
{\bf ab} and the {\bf cd}-index
and explain why one expects results
about the {\ab}-index to carry over to the {\bf cd}-index.
In Appendix~\ref{s:recur},
we give a recursion for computing
the $\cd$-index of an Eulerian poset
in terms of certain polynomial sequences.
These may be of independent interest.
Appendix~\ref{s:cube} deals with the cubical lattice,
which is the face lattice of the cube.
Just as the simplex is the simplest polytope,
the cube is the simplest zonotope
and is an object of interest in its own right.
Usually, techniques that 
work for the Boolean lattice also work for the 
cubical lattice with minor modifications; see
~\cite{\emahajan,\ereaddy}.
Following this general principle, 
we establish similar results for the {\bf cd}-index 
of the cubical lattice. 
In Appendix~\ref{s:fun}, 
we show that $\Fhat$ is a free algebra 
with the $\centerdot$ product (Theorem~\ref{t:free}).
We also show that the Dehn-Sommerville relations 
satisfied by the flag $f$-vector of an Eulerian poset
are equivalent to certain  simple identities 
that hold in $\Fhat$.

\section{The polynomial algebra $k \langle \cv,\dv \rangle$}
\label{s:algebra}

The basic algebraic object to consider is
${\mathcal A} = k \langle \av,\bv \rangle$, 
the free algebra in two non-commuting variables {\bf a} and {\bf b}.
The other object of interest is the subalgebra 
${\mathcal F}$,
generated by the elements 
${\bf c} := \av + \bv$ and 
${\bf d} := \ab + \bv\av$. 
Since we are primarily interested in the $\cd$-index,
we will concentrate on $\Ff$ and never deal with $\mathcal A$.

\subsection{The basic setup}

We begin by recalling some facts from~\cite{\ereaddy}.
Let $k$ be a field of characteristic 0.
Let ${\mathcal F} = k \langle \cv,\dv \rangle$ 
be the polynomial algebra in the non-commuting
variables {\bf c} and {\bf d}. 
By setting the degree of {\bf c} to be
1 and of {\bf d} to be $2$, we write
${\mathcal F} = \oplus_{n \geq 0} {\mathcal F}_n$, 
where ${\mathcal F}_n$ is spanned by the {\bf cd}-monomials of degree $n$.
The product in $\Ff$ is given by concatenation 
and the unit element is 1.

\begin{proposition}[Ehrenborg-Readdy] \label{p:d}
The vector space ${\mathcal F}$
has a (coassociative) coproduct
$\Delta \colon \mathcal F \rightarrow \mathcal F \tensor \mathcal F$ given by the
initial conditions
$$\Delta(1) = 0, \: \: \Delta({\bf c}) = 2  (1 \tensor 1), 
\: \: \Delta({\bf d}) = 1 \tensor {\bf c} + {\bf c} \tensor 1$$
and the rule
$\Delta(u  v) = \Delta(u)  v + u  \Delta(v)$
for $u,v \in \Ff.$
\end{proposition}

Under the coproduct $\D$, 
the vector space $\Ff$ is a coassociative coalgebra,
but without a counit map.
Further, the rule says that $\D$ is a derivation on $\Ff$
into the
$(\Ff,\Ff)$-bimodule $\mathcal F \tensor \mathcal F$.
This makes $\Ff$ an infinitesimal bialgebra,
also called a Newtonian coalgebra.
This notion was first defined by Joni and Rota~\cite{\joni}.
For more recent work, 
see the papers of Aguiar~\cite{\marcelocd,\marcelo}. 
However, we will not
use any facts about infinitesimal bialgebras.

The motivation for the definition of $\D$
is as follows.
Consider the map 
\[
\Psi : \{ \text{Eulerian posets} \} \rightarrow k \langle \cv , \dv \rangle,
\]
which assigns to an Eulerian poset $P$ its $\cd$-index $\Psi(P)$.
The vector space spanned by all Eulerian posets
is a coalgebra with the coproduct given by
\[
\Delta(P) = \sum_{\hz < x < \ho} [\hz,x] \tensor [x,\ho].
\]
And the map $\Psi$ is a morphism of coalgebras.
In other words,
the identity
\begin{equation} \label{e:coalg_morph}
\Delta(\Psi(P)) = \sum_{\hz < x < \ho} \Psi([\hz,x]) \tensor \Psi([x,\ho])
\end{equation}
holds for any Eulerian poset $P$;
see~\cite[Proposition 3.1]{\ereaddy}.

We will use this later to derive a basic result
about the Boolean lattice; see Lemma~\ref{l:easy}.
We will also return to it briefly in Appendix~\ref{s:fun},
when we discuss the Dehn-Sommerville equations.

\begin{proposition}[Ehrenborg-Readdy] \label{p:g}
There is a well-defined linear map $G \colon \Ff \rightarrow \Ff$
given by the initial conditions
$$G(1) = 0, \: \: G(\cv) = \dv, \: \: G(\dv) = \cv \dv$$
and the rule
$G(u v) = G(u) v + u G(v)$, 
such that 
$$\Psi(B_{n+1}) = \Psi(B_n) \cv + G(\Psi(B_n)).$$
\end{proposition}
\noindent
The importance of the map $G$ is that 
it gives an inductive way of computing $\Psi(B_n)$.

\subsection{An extension of the basic setup} \label{ss:ext}

Consider
$\Fhat = k  \ev \oplus \Ff$, where
$\ev$ is a formal symbol with degree $-1$. 
We write
$\Fhat = \oplus_{n \geq -1} {\mathcal F}_n$, 
where $\Ff_{-1} = k  \ev$.
Define a coproduct
$\Dhat \colon \Fhat \rightarrow \Fhat \tensor \Fhat$
by
\begin{equation}
\Dhat(\ev) = \ev \tensor \ev \quad \text{and} \quad
\Dhat(u) = \D(u) + \ev \tensor u + u \tensor \ev,
\end{equation}
for $u \in \Ff$. 
Observe that
$\Dhat$ has degree $-1$, that is,
$\Dhat({\mathcal F}_{n}) \subseteq 
\bigoplus_{i \geq -1} {\mathcal F}_i \tensor {\mathcal F}_{n-i-1}.$
Also let
$\varepsilon \colon \Fhat \rightarrow k$ be given by 
the delta function $\delta_\ev$.
It is easy to see that $\Fhat$ is a coalgebra, with $\Dhat$ as the
coproduct and $\varepsilon$ as the counit. 
The process of passing from $\Ff$ to $\Fhat$
just described is the standard way of
adding a counit to a coalgebra.
For matters of notational convenience,
we let
$\ev v = v \ev = 0.$
With this convention, 
it is still true that
$\Dhat(u v) = \Dhat(u) v + u \Dhat(v)$
for $u,v \in \Fhat$.
Hence the extended object $\Fhat$ 
is also an infinitesimal bialgebra.

\begin{figure} 
$$  \begin{array}{c | c c c}
      & \Dhat & \hat G & \hat H \\ \hline
\ev   & \ev \tensor \ev & 1 & \\
1     & \ev \tensor 1 + 1 \tensor \ev & \cv & \cv \\
\cv   & 2 (1 \tensor 1) + \cv \tensor \ev + \ev \tensor \cv  
                                           & \cv^2 + \dv & \cv^2 + 2 \dv \\
\dv   & 1 \tensor \cv + \cv \tensor 1 + \dv \tensor \ev + \ev \tensor \dv  
          & \cv \dv + \dv \cv & \dv + 2 \dv \cv 
     \end{array} $$
\caption{Values of the maps $\Dhat,\hat G$ and $\hat H$
          at $\ev,1,\cv,\dv$.}
\label{f:val}
\end{figure}

Let $\hat G \colon \Fhat \rightarrow \Fhat$
be the linear map defined by ${\hat G}(\ev) = 1$
and 
${\hat G}(u) = G(u) + u \cv$ for $u \in \Ff$.
Note that the definition of $\hat G$ 
is arranged so that the equation
\begin{equation} \label{e:g}
{\hat G}(\Psi(B_n)) = \Psi(B_{n+1}) \quad \text{holds for} \quad n \geq 0.
\end{equation}
This can be seen from Proposition~\ref{p:g}.
It is clear that understanding $\hat G$ is crucial for our
purposes. 
At least, that is the approach we take.

It is easy to see that 
for $u,v \in \Ff$, we have
$\hat G(u v) = \hat G(u) v + u \hat G(v) 
        - u \cv v$.
We will use this identity later in the proof of
Theorem~\ref{t:coderivation}. 
However, it does not hold in $\Fhat.$ 
Hence the terms 
$\hat G(\ev v)$ and
$\hat G(u \ev)$ need to be handled with care.

Figure~\ref{f:val} illustrates the maps 
$\Dhat$ and $\hat G$.
The map $\hat H$ is the analogue of the map $\hat G$
for the cubical lattice.
This will be explained in Appendix~\ref{s:cube}.

\subsection{More definitions}

Apart from $\Dhat$ and $\hat G$, 
we define a third map
$\mu \colon \Fhat \tensor \Fhat \rightarrow \Fhat$ 
of degree $2$ by
$\mu(\ev \tensor \ev) = 2, \: \mu(\ev \tensor v) = \cv v, \:
\mu(v \tensor \ev) = v \cv$ and 
$\mu(u \tensor v) = u \dv v$ for $u,v \in \Ff$.
The relation of $\mu$ with the previous two maps
is given by Lemma~\ref{l:coderivation_identity}.

For $v$, a {\bf cd}-monomial of degree $n$,
let $\beta(v)$ be the coefficient of 
$v$ in $\Psi(B_{n+1})$.
In more fancy language,
$\beta(v) = \langle \delta_v,\Psi(B_{n+1}) \rangle$.
Observe that $\beta(\ev) = 1$.
We then extend the definition to $\Fhat$ by linearity.

Let $v^*$ denote the reverse of the {\bf cd}-monomial $v$.
Also define $(u \tensor v)^*$ to be
$v^* \tensor u^*$.

It is easy to check from the definitions that
$\Dhat(u^*) = (\Dhat(u))^*$,
$\hat G(u^*) = (\hat G(u))^*$ and
$\mu((u \tensor v)^*) = (\mu(u \tensor v))^*$.
It is known that 
$\beta(v) = \beta(v^*)$.
This also follows by induction from the second equality
and equation~\eqref{e:g}.

\subsection{The main result}

We now state and prove the main result of this section.

\begin{theorem} \label{t:coderivation}
Let $\Dhat$ and $\hat G$ be as defined before. Then
\[
\Dhat \circ {\hat G} = (\id \tensor {\hat G} + {\hat G} \tensor \id)
\circ \Dhat.
\]
\noindent
In other words,
${\hat G}$ is a coderivation on $\Fhat$
with respect to $\Dhat$.
\end{theorem}
\begin{proof}
We evaluate both sides of the identity at 
an arbitrary {\bf cd}-monomial and then use
induction on its degree to show that 
they yield the same value.

The first step is to check directly that the 
identity holds at the {\bf cd}-monomials
$\ev, 1, \cv$ and $\dv$.
This is straightforward to check using
Figure~\ref{f:val}.
To complete the induction step, 
we begin by evaluating the RHS at $u v$ and 
expanding using the inductive definitions of 
$\Dhat$ and $\hat G$. The induction step is shown below.
\[
\begin{array}{r c l}
RHS_{u v} & = & (\id \tensor {\hat G} + {\hat G} \tensor \id)
(\Dhat(u v)) \\
                & = & (\id \tensor {\hat G} + {\hat G} \tensor \id)
(\Dhat(u) v + u \Dhat v).
\end{array}
\]
On expanding further, we obtain four terms, 
two of which we write down explicitly.
\[
\begin{array}{r c l}
(\id \tensor {\hat G})(\Dhat(u) v) 
    & = & (\id \tensor {\hat G})(\Dhat(u)) v +
\Dhat(u) \hat G(v) -
\Dhat(u) \cv v -
u \tensor v,\\
({\hat G} \tensor \id)(\Dhat(u) v)
    & = & ({\hat G} \tensor \id)(\Dhat(u)) v.
\end{array}
\]
The correction term
$- u \tensor v$ 
in the first expression 
accounts for the difference between the terms
$(\id \tensor {\hat G})(u \tensor \ev v)$ and
$(\id \tensor {\hat G})(u \tensor \ev) v.$  

The remaining two terms can be written down by symmetry.
Summing up all the four terms and
applying induction, we obtain
\[
\begin{array}{r c l}
RHS_{u v} 
    & = & \Dhat(\hat G(u)) v + u \Dhat(\hat G(v)) + 
          \hat G(u) \Dhat(v) + \Dhat(u) \hat G(v) \\ 
    &   & \: \: \: - \Dhat(u) \cv v - u \cv \Dhat(v) -
           2 u \tensor v\\
    & = & \Dhat(\hat G(u) v) + \Dhat(u \hat G(v)) -
           \Dhat(u \cv v) \\
    & = & \Dhat(\hat G(u v)) \\
    & = & LHS_{u v}.
\end{array}
\]
\end{proof}

\begin{lemma} \label{l:coderivation_identity}
Let $\hat G$, $\mu$ and $\Dhat$ be as defined before.
Then
$2 \hat G = \mu \circ \Dhat.$
\end{lemma}
\begin{proof}
The proof follows the same pattern 
as that of the previous theorem.
The induction step is as follows.
\[
\begin{array}{r c l}
RHS_{u v} 
    & = & \mu(\Dhat(u) v + u \Dhat (v)) \\
    & = & \mu(\Dhat(u)) v + u \mu(\Dhat (v)) -
          2 u \cv v \\
    & = & 2 (\hat G(u) v + u \hat G(v) 
        - u \cv v) \\
    & = & 2 (\hat G(u v)) \\
    & = & LHS_{u v}.
\end{array}
\]
\end{proof}

\begin{remark}
One may check that the map $\hat G$ 
is also a derivation with respect to the product $\mu$.
And the triple $(\Fhat, \Dhat, \mu)$
is an infinitesimal bialgebra.
This was pointed out by Marcelo Aguiar.
\end{remark}

\section{The dual setup} \label{s:duality}

In this section,
we present the picture dual to the one in Section~\ref{s:algebra}.
For motivation, we give some simple applications 
in Section~\ref{ss:sa}.
In Section~\ref{ss:explicit},
we write down explicit formulas for the dual maps.
These lead to some immediate consequences,
which we discuss in Section~\ref{ss:ssa}.

\subsection{The product $\centerdot$ and the derivation $S$}
\label{ss:duality}
Let 
${\Fhat^{*}}$
be the restricted dual of
${\Fhat}$,
namely the space of linear functionals on $\Fhat$
that vanish 
on the graded piece $\Ff_n$ for sufficiently large $n$. 
As noted before, ${\Fhat}$ has a basis consisting 
of all the {\bf cd}-monomials. This gives 
${\Fhat^{*}}$ 
a natural basis consisting of the delta functions 
$\delta_v$, where $v$ is any {\bf cd}-monomial.
As vector spaces,
$\Fhat$
and
$\Fhat^*$
are isomorphic and we identify them
using our specific choice of bases
$v \leftrightarrow \delta_v$.
This may look unnatural at first, 
but it is not so strange given that 
we are indeed biased towards a particular basis for ${\Fhat}$ and 
are trying to study the {\bf cd}-index in this basis.

Dualise the maps 
$\Dhat, \hat G$ and $\mu$ 
defined in the previous section to get
the corresponding dual maps
$\Dhat^*, \hat G^*$ and $\mu^*$.
The map
$\Dhat^*$
is the convolution product on 
$\Fhat^*$,
which was first introduced by Kalai~\cite{\kalai}.
The algebra $(\Fhat^*,\Dhat^*)$
can be identified with the algebra $A_{\varepsilon}$
studied by Billera and Liu~\cite{\bn}.

Now using the identification of 
$\Fhat^*$
with
$\Fhat$ explained above,
we transfer these maps back to $\Fhat$ and
obtain three maps, 
\[
\centerdot \colon \Fhat \tensor \Fhat \rightarrow \Fhat, \
S \colon \Fhat \rightarrow \Fhat \ \text{and} \
\mu^* \colon \Fhat \rightarrow \Fhat \tensor \Fhat.
\]
\vanish{
$\centerdot \colon \Fhat \tensor \Fhat \rightarrow \Fhat$, \
$S \colon \Fhat \rightarrow \Fhat$ and
$\mu^* \colon \Fhat \rightarrow \Fhat \tensor \Fhat.$
}
We note that
these maps have degrees $1$, $-1$ and $-2$ respectively.
In other words,
\[
\centerdot \colon \Ff_i \tensor \Ff_j \rightarrow \Ff_{i+j+1}, \
S \colon \Ff_{i+1} \rightarrow \Ff_i \ \text{and} \
\mu^* \colon \Ff_n \rightarrow \bigoplus\limits_{i \geq -1} \Ff_i \tensor \Ff_{n-i-2}.
\]
\vanish{
$\centerdot \colon \Ff_i \tensor \Ff_j \rightarrow \Ff_{i+j+1}$, \
$S \colon \Ff_{i+1} \rightarrow \Ff_i$ and
$\mu^* \colon \Ff_n \rightarrow \bigoplus\limits_{i \geq -1} \Ff_i \tensor \Ff_{n-i-2}.$
}
The definitions of these maps can be made very explicit,
see Section~\ref{ss:explicit}.
But before doing that, 
we will present the dual versions of the results of the
previous section and derive some immediate consequences from them.
This would give some motivation for considering these
dual maps.

By general principles of duality,
$\Fhat$ is an associative algebra with unit $\ev$, 
with respect to the $\centerdot$ product.
The dual versions of
Theorem~\ref{t:coderivation} 
and Lemma~\ref{l:coderivation_identity} 
are as follows.

\begin{theorem} \label{t:derivation}
The map $S$ is a derivation on the algebra $\Fhat$, that is,
\[
S \circ \centerdot = \centerdot \circ (\id \tensor S + S \tensor \id).
\]
This may be more familiarly
expressed as 
$S(u \centerdot v) = S(u) \centerdot v + u \centerdot S(v)$
for $u,v \in \Fhat$.
Also,
$S(1) = \ev$ and 
$S(\ev) = 0$.
\end{theorem}

\begin{lemma} \label{l:derivation_identity}
Let the maps $S$, $\centerdot$ and $\mu^*$ be as defined before. Then 
$2 S = \centerdot \circ \mu^*$.
\end{lemma}

We know from equation~\eqref{e:g} that
$\hat G$ satisfies the important property
${\hat G}(\Psi(B_n)) = \Psi(B_{n+1})$.
We now state the dual version of this property.

\begin{lemma} \label{l:reduction}
Let $v$ be any {\bf cd}-monomial of non-negative degree. Then
$\beta(S(v)) = \beta(v)$.
\end{lemma}
\begin{proof}
To illustrate how duality works,
we give a proof of this lemma. 
Let $v$ be a monomial of degree $n$, with $n \geq 0$. Then,
$$  \begin{array}{r c c c l}
\beta(S(v)) & = & \langle \delta_{S(v)},\Psi(B_{n}) \rangle                          
     & = &
           \langle {\hat G}^*(\delta_v),\Psi(B_n) \rangle \\           
    & & 
    & = &
           \langle \delta_v,{\hat G}(\Psi(B_n)) \rangle \\  
    & &
     & = &
           \langle \delta_v,\Psi(B_{n+1}) \rangle \\
    & &
     & = &
           \beta(v),
     \end{array} $$
where the second last equality uses the identity 
${\hat G}(\Psi(B_n)) = \Psi(B_{n+1})$.
\end{proof}

\subsection{Simple applications} \label{ss:sa}

We now show some interesting consequences of the ideas
discussed so far.

\begin{lemma} \label{l:easy}
Let $u$ and $v$ be $\cd$-monomials of degree $m$ and $n$ 
respectively. Then
\[
\beta(u \centerdot v) = \binom {m+n+2}{m+1} \beta(u) \beta(v).
\]
\end{lemma}
\begin{proof}
The key fact to use is that
an interval in a Boolean lattice
is again a smaller Boolean lattice.
The dual to equation~\eqref{e:coalg_morph} is the identity
$$
\Psi_P^*(u \centerdot v) = \sum_{\hz < x < \ho} 
\Psi_{([\hz,x])}^*(u)  \Psi_{([x,\ho])}^*(v),
$$
where $\Psi_{P}^*(w)$
denotes the coefficient of $w$ in $\Psi(P)$.
Setting $P$ to be the Boolean lattice $B_{m+n+2}$,
we obtain
$$
\beta(u \centerdot v) = \sum_{\hz < x < \ho} \beta(u) \beta(v),
$$
where the sum is over those $x$'s in $B_{m+n+2}$
whose rank is $m+1$.
The result now follows.
\end{proof}

\begin{example} \label{eg:warmup}
To illustrate
how Theorem~\ref{t:derivation}
and Lemma~\ref{l:reduction} 
work together,
we compute
$\beta(\underbrace {1 \centerdot 1 \centerdot \ldots \centerdot 1}_{n}).$
First by Lemma~\ref{l:reduction}, we have
$\beta(\underbrace {1 \centerdot 1 \centerdot \ldots \centerdot 1}_{n}) =
\beta(S(\underbrace {1 \centerdot 1 \centerdot \ldots \centerdot 1}_{n})).$
Next by Theorem~\ref{t:derivation},
$S(\underbrace {1 \centerdot 1 \centerdot \ldots \centerdot 1}_{n}) = 
n S(\underbrace {1 \centerdot 1 \centerdot \ldots \centerdot 1}_{n-1}).$
These two facts
and the initial condition 
$\beta(1) = 1$, 
yield us the result 
$\beta(\underbrace {1 \centerdot 1 \centerdot \ldots \centerdot 1}_{n}) = n!.$

\end{example}

Next we record two results
that will be useful in later sections.
They are direct corollaries of Lemma~\ref{l:easy}
and the fact $\beta(v) = \beta(v^*)$.
However, we will not rely on this lemma.
Instead, we will give independent proofs 
using the method illustrated in Example~\ref{eg:warmup}.
For that, observe two simple facts, namely,
$S(u)^* = S(u^*)$ and 
$(u \centerdot v)^* = v^* \centerdot u^*$.
Also recall that the map $S$ has degree $-1$.

\begin{lemma} \label{l:add_on}
Let $u$ and $v$ be {\bf cd}-monomials of the same degree
and $w$ be any {\bf cd}-monomial. Then we have,
$$
\begin{array}{r c l}
\beta(u) > \beta(v) & \mbox{\rm iff} 
         & \beta(u \centerdot w) > \beta(v \centerdot w) \\
\beta(u) = \beta(v) & \mbox{\rm iff} 
         & \beta(u \centerdot w) = \beta(v \centerdot w).
\end{array}
$$
\end{lemma}
\begin{proof}
We first prove the forward implications of both statements.
Since the proofs are similar, we only do the forward 
implication of the first statement.

Perform induction on the degrees of $u$,$v$ and $w$.
The induction base is straightforward.
Now consider the induction step.
By our assumption,
$\beta(u) > \beta(v)$
and hence
$\beta(S(u)) > \beta(S(v))$.
Therefore by induction we obtain,
$\beta(u \centerdot S(w)) > \beta(v \centerdot S(w))$
and
$\beta(S(u) \centerdot w) > \beta(S(v) \centerdot w)$.
Summing up and using Theorem~\ref{t:derivation} and
Lemma~\ref{l:reduction}, we get
$\beta(u \centerdot w) > \beta(v \centerdot w).$

The backward implications of both the statements are
again similar. To see the backward implication of
the first statement, 
assume the contrary,
that is, either
$\beta(u) < \beta(v)$ or
$\beta(u) = \beta(v)$. 
Then the forward implications, which we 
just  proved, give a contradiction.
\end{proof}

\begin{lemma} \label{l:switch}
Let $u$ and $v$ be {\bf cd}-monomials. Then we have,
$\beta(u \centerdot v) = \beta(u \centerdot v^*) = \beta(u^* \centerdot v)$
and
$\beta(u \centerdot v) = \beta(v \centerdot u)$.
\end{lemma}
\begin{proof}
We first prove the first statement by an induction on 
the size of the {\bf cd}-monomial. It is enough to show only
the first equality. The second one follows by the symmetry
in our argument.
The induction base is provided by the statements
$\beta(\ev \centerdot v) = \beta(\ev \centerdot v^*)$ and
$\beta(u \centerdot \ev) = \beta(u \centerdot \ev^*)$.
$$  
\begin{array}{r c l}
\beta(u \centerdot v) 
    & = &
           \beta(S(u \centerdot v)) \\
    & = &
           \beta(S(u) \centerdot v) + \beta(u \centerdot S(v)) \\
    & = &
           \beta(S(u) \centerdot v^*) + \beta(u \centerdot S(v^*)) \\
    & = &
           \beta(S(u \centerdot v^*)) \\
    & = &
           \beta(u \centerdot v^*).
\end{array} 
$$
We made use of the induction hypothesis in the third step.

The second statement follows from the first
by the chain of inequalities shown below.
$$ 
\beta(u \centerdot v) = \beta((u \centerdot v)^*) = 
        \beta(v^* \centerdot u^*) = \beta(v \centerdot u).
$$
\vanish{
$$  \begin{array}{r c c c c c l}
\beta(u \centerdot v) 
    & = &
           \beta((u \centerdot v)^*) 
    & = &
           \beta(v^* \centerdot u^*) 
    & = &
           \beta(v \centerdot u).
     \end{array} $$
}
\end{proof}

\subsection{An explicit description of the maps $\centerdot, S$ and $\mu^*$}
\label{ss:explicit}

We obtained the maps $\centerdot, S$ and $\mu^*$ 
by taking duals of certain other maps. 
In this section,
we go through the duality grind to
give explicit formulas for these maps.
By way of justification, we give some
applications in Section~\ref{ss:ssa}.

\begin{lemma} \label{l:bullet}
The $\centerdot$ product on $\Fhat$ is determined by the initial conditions
\[
1 \centerdot 1 = 2 {\bf c}, \: \:
1 \centerdot {\bf c} = {\bf c} \centerdot 1 = {\bf d} + 2 {\bf c}^{2}, \: \:
1 \centerdot {\bf d} = 2 {\bf c} {\bf d}, \: \:
{\bf d} \centerdot 1 = 2 {\bf d} {\bf c},
\]
\[
{\bf c} \centerdot {\bf c} = {\bf d} {\bf c} + {\bf c} {\bf
d} + 2 {\bf c}^{3}, \: \:
{\bf c} \centerdot {\bf d} = {\bf d}^{2} + 2 {\bf c}^{2} {\bf d}, \: \:
{\bf d} \centerdot {\bf c} = {\bf d}^{2} + 2 {\bf d} {\bf
c}^{2}, \: \:
{\bf d} \centerdot {\bf d} = 2 {\bf d} {\bf c} {\bf d},
\]

\medskip
\noindent 
and the rule
$(u \epsilon_1) \centerdot (\epsilon_2 v) =
u (\epsilon_1 \centerdot \epsilon_2) v$,
where 
$\epsilon_1$ 
and 
$\epsilon_2$ 
are either of the letters $\cv,\dv$ and $u,v$ are {\bf cd}-monomials.
\end{lemma}
\begin{proof}
To show 
$\cv \centerdot \dv = 2 \cv^2 \dv + \dv^2$, 
for example,
we prove the equivalent statement,
$\delta_{\cv \centerdot \dv} =
2 \delta_{\cv^2 d} + \delta_{\dv^2}$.
Evaluating the LHS at the {\bf cd}-monomial $w$, we obtain
\[
\langle\delta_{\cv \centerdot \dv},w\rangle =
\langle\Dhat^*(\delta_{\cv} \tensor \delta_{\dv}),w\rangle = 
\langle\delta_{\cv} \tensor \delta_{\dv},\Dhat(w)\rangle =  
\langle2 \delta_{\cv^2 d} + \delta_{\dv^2},w\rangle.
\]

\vanish{
$$  \begin{array}{r c c c c c c}
\langle\delta_{\cv \centerdot \dv},w\rangle 
     & = & \langle\Dhat(\delta_{\cv} \tensor \delta_{\dv}),w\rangle   
     & = & \langle\delta_{\cv} \tensor \delta_{\dv},\Dhat(w)\rangle   
     & = & \langle2 \delta_{\cv^2 d} + \delta_{\dv^2},w\rangle.
     \end{array} $$
}

\noindent
The last equality is true since 
$\cv^2 \dv$ and $\dv^2$ 
are the only monomials whose coproduct involves the term
$\cv \tensor \dv$.
The other verifications are similar and the reader may try out a few
to get a feel for this product.

To check the rule stated in the lemma, we show
$\delta_{(u \epsilon_1) \centerdot (\epsilon_2 v)} =
\delta_{u (\epsilon_1 \centerdot \epsilon_2) v}$.
To do this, evaluate the LHS at the {\bf cd}-monomial $w$. 
Also assume that $w = u w' v$ for some 
{\bf cd}-monomial $w'$.
If $w$ does not have this form, then both sides evaluate to zero.
$$
\begin{array}{r c l}
\langle \delta_{(u \epsilon_1) \centerdot (\epsilon_2 v)}, \: w \rangle
    & = &
          \langle \Dhat^*(\delta_{(u \epsilon_1)} \tensor
\delta_{(\epsilon_2 v)}), \: w \rangle \\
    & = &
          \langle \delta_{(u \epsilon_1)} \tensor
\delta_{(\epsilon_2 v)}, \: \Dhat(w) \rangle \\
    & = &
          \langle \delta_{\epsilon_1} \tensor
\delta_{\epsilon_2}, \: \Dhat(w') \rangle \\
    & = &
          \langle \delta_{\epsilon_1 \centerdot
\epsilon_2}, \: w' \rangle \\
    & = &
          \langle \delta_{u (\epsilon_1 \centerdot
\epsilon_2) v}, \: w \rangle.
\end{array}
$$
The third equality follows from the rule for $\Delta$
stated in Proposition~\ref{p:d}.
\end{proof}
\noindent
Next we describe the maps
$S$ and $\mu^*$.
The proofs are straightforward and we omit them.
Note that every {\bf cd}-monomial
can be uniquely written in the form
$\cv^{m_1} \dv \cv^{m_2} \dv \ldots \dv \cv^{m_k},$
where 
$m_1,m_2,\ldots,m_k$
are non-negative integers.

\begin{lemma} \label{l:definition_S}
Let $m_1,m_2,\ldots,m_k$
be non-negative integers. The map $S$ is given by
$
S(\cv^{m_1} \dv \cv^{m_2}  \dv \ldots \dv  \cv^{m_k}) =
\sum_{i=1}^k \cv^{m_1} \ldots \dv  \cv^{m_i-1}  \dv \ldots \cv^{m_k} +\\
\hspace {20 mm} \sum_{i=1}^{k-1} \cv^{m_1} \ldots \dv  \cv^{m_i}  \cv  \cv^{m_{i+1}} \dv \ldots \cv^{m_k}.
$
\end{lemma}

\begin{lemma}
Let $m_1,m_2,\ldots,m_k$
be non-negative integers. The map $\mu^*$ is given by
$
\mu^*(\cv^{m_1} \dv \cv^{m_2}  \dv \ldots \dv  \cv^{m_k}) =
\ev \tensor (\cv^{m_1-1} \dv \ldots \dv  \cv^{m_k}) +
(\cv^{m_1} \dv \ldots \dv  \cv^{m_k-1}) \tensor \ev + \\
\sum_{i=1}^{k-1} \cv^{m_1} \ldots \dv  \cv^{m_i} 
\tensor \cv^{m_{i+1}} \dv \ldots \cv^{m_k}.
$

\vanish{
$
\begin{array}{r c l}
\mu^*(\cv^{m_1} \dv \cv^{m_2}  \dv \ldots \dv  \cv^{m_k}) 
   & = &
\ev \tensor (\cv^{m_1-1} \dv \ldots \dv  \cv^{m_k}) +
(\cv^{m_1} \dv \ldots \dv  \cv^{m_k-1}) \tensor \ev + \\
   &   & 
\sum_{i=1}^{k-1} \cv^{m_1} \ldots \dv  \cv^{m_i} 
\tensor \cv^{m_{i+1}} \dv \ldots \cv^{m_k}.
\end{array}
$
}
\end{lemma}
\noindent
Combining this lemma with 
Lemma~\ref{l:derivation_identity},
we obtain a more useful expression for $S$
as follows.

\begin{lemma} \label{l:formula_useful}
The map $S$ is given by the equation,
$2 S(\cv^{m_1} \dv \cv^{m_2}  \dv \ldots \dv  \cv^{m_k}) =
(\cv^{m_1-1} \dv \ldots \dv  \cv^{m_k}) +
(\cv^{m_1} \dv \ldots \dv  \cv^{m_k-1}) + 
S'(\cv^{m_1} \dv \cv^{m_2}  \dv \ldots \dv  \cv^{m_k})$,
where \\
$S'(\cv^{m_1} \dv \cv^{m_2}  \dv \ldots \dv  \cv^{m_k}) = 
\sum_{i=1}^{k-1} \cv^{m_1} \ldots \dv  \cv^{m_i} \centerdot \cv^{m_{i+1}} \dv 
\ldots \cv^{m_k}.$
\end{lemma}
\noindent
The above lemma can also be checked directly from 
Lemmas~\ref{l:bullet} and~\ref{l:definition_S}.

\subsection{More applications} \label{ss:ssa}

The explicit formula given by Lemma~\ref{l:bullet}
allows us to make more concrete sense out of
Lemma~\ref{l:switch}.

\begin{lemma} \label{l:identities}
Let $u$, $v$ and $w$ be any {\bf cd}-monomials. Then we have,
\[
\beta(\dv u \dv \cv \dv v) =
\beta(\dv u^* \dv \cv \dv v) \quad \text{and} \quad
\beta(u \dv \cv \dv 
v \dv \cv \dv w) =
\beta(u \dv \cv \dv 
v^* \dv \cv \dv w).
\]
\end{lemma}
\begin{proof}
By Lemma~\ref{l:bullet},
observe that
$(\dv u^* \dv \cv \dv v) =
1/2 (\dv u^* \dv) \centerdot (\dv v)$.
Now, the first identity follows from the following
sequence of equalities.
\[
\beta(\dv u \dv \cv \dv v) =
1/2 \beta((\dv u \dv) \centerdot (\dv v)) =
1/2 \beta((\dv u^* \dv) \centerdot (\dv v)) =
\beta(\dv u^* \dv \cv \dv v).
\]
The second inequality follows from 
Lemma~\ref{l:switch}.

The second result can be proved similarly from the identity
$(u \dv \cv \dv 
v \dv \cv \dv w) =
    1/4 (u \dv) \centerdot (\dv v \dv) 
            \centerdot (\dv w).$
\end{proof}

These identities look exciting and one may ask for a complete list
of such identities. We do not attempt to answer this question.
As an interesting fact,
direct computation shows that for Boolean lattices of rank at most 13,
there is only one identity that 
Lemma~\ref{l:identities}
does not account for,
namely,
$\beta(\cv^2 \dv^2 \cv^3 \dv \cv) =
\beta(\cv^3 \dv \cv^2 \dv \cv \dv).$
We have no explanation for this identity
or any others of this type that might exist.
A similar but more complicated argument leads to 
a different class of identities
that we will do in Corollary~\ref{c:identity}.

Now we provide some examples of how
Lemmas~\ref{l:reduction}
and
~\ref{l:definition_S}
work together.

\begin{example} \label{eg:compute}
Using the reduction
$S(\cv^m) = \cv^{m-1}$
and the fact that $\beta(1)=1$,
we obtain the known
fact that 
$\beta(\cv^m) = 1$ for all $m$.
Next we compute
$\beta(\cv^i \dv \cv^j).$
Using
Lemmas~\ref{l:reduction}
and
~\ref{l:definition_S},
we get
$\beta(\cv^i \dv \cv^j) = \beta(\cv^{i-1} \dv \cv^j) +
\beta(\cv^i \dv \cv^{j-1}) +1.$
We rewrite this equation as
$(\beta(\cv^i \dv \cv^j) + 1) = 
(\beta(\cv^{i-1} \dv \cv^j) + 1) +
(\beta(\cv^i \dv \cv^{j-1}) +1)$,
which reminds us of the binomial recursion
${\binom {n}{k}} = {\binom {n-1}{k-1}} + {\binom {n-1}{k}}.$
After checking the initial conditions, we conclude that
$\beta(\cv^i \dv \cv^j) = {\binom {i+j+2}{i+1}} -1.$
\end{example}

\begin{lemma} \label{l:cc_d}
Let $u$ and $v$ be {\bf cd}-monomials. Then
$\beta(u \dv v) \geq \beta(u \cv^2 v)$, 
with equality if $u$ and $v$ are both empty.
\end{lemma}
\begin{proof}
We do an induction on the degree of the
{\bf cd}-monomial. 
By Lemma~\ref{l:reduction}, 
it is enough to show the equivalent statement
$\beta(S(u \dv v)) \geq \beta(S(u \cv^2 v))$.
Write 
$u = u' \cv^m$ and
$v = \cv^n v'$,
where $u'$ ends with a $\dv$ 
and $v'$ begins with a $\dv.$
Now observe that 
$S(u \dv v) = S(u') \cv^m \dv \cv^n v' +
u' S(\cv^m \dv \cv^n) v' + u' \cv^m \dv \cv^n S(v').$
A similar expansion can be written out for
$S(u \cv^2 v).$
Applying the induction hypothesis on the first and third terms
and using Lemma~\ref{l:definition_S} on the second term,
the desired result follows.
\end{proof}

The reader will notice a common method in our examples.
In order to prove any result about $\beta(v)$, we start by
looking at $\beta(S(v))$. Then we expand $S(v)$ directly
using the description given by Lemma~\ref{l:definition_S}.
Since $S$ has degree $-1$, every term that occurs 
has degree lower than $v$.
We then group terms together such that
every grouped term satisfies the induction hypothesis.
And the result gets proved by induction.

The procedure above is a brute force technique and can
be messy on more complicated examples.
Hence we prefer to use the description of $S$ provided by 
Lemma~\ref{l:formula_useful}.
This involves the $\centerdot$ product and
therefore gives us access to Lemmas~\ref{l:add_on}
and~\ref{l:switch}.
There are many interesting inequalities which can be derived
from Lemma~\ref{l:formula_useful}. However, it is easier to express 
them in an alternate notation,
which we define in the next section.

\section{The list notation} \label{s:ln}

As may have become evident by now,
it is easier to work with a more compact notation,
where we use lists to denote {\bf cd}-monomials.
Apart from the advantage of being compact,
it is well suited for the definition of
the maps
$\centerdot, S$ and $\mu^*.$
In Section~\ref{ss:list},
we develop this notation and
then restate all the important results in terms of lists.
In Section~\ref{ss:ma},
we give some further results.

\subsection{Restatement of results using the list notation}
\label{ss:list}

Every {\bf cd}-monomial
can be uniquely written in the form
$\cv^{m_1} \dv \cv^{m_2} \dv \ldots \dv \cv^{m_k},$
where 
$m_1,m_2,\ldots,m_k$
are non-negative integers.
Hence we may represent it by the list
$(m_1,m_2,\ldots,m_k).$
We define the length of a list to be the number of elements in it.
Note that the list
$(m_1,m_2,\ldots,m_k)$
has degree
$(\sum_{i=1}^k m_i) + 2 (k-1)$
and length $k$.
The element 1 is denoted by the list $(0)$ 
and the unit element $\ev$ by the empty list.
We follow the convention that a list cannot
have negative entries. 
If a list with negative entries appears in a definition
or computation, then we simply ignore it.
In other words, we define it to be zero.
Also for a list $M$, let $M^*$ denote its reverse.
For future convenience, we now restate
the results of the last section
using lists.

\begin{lemma} \label{l:add_on_list}
Let $L$ and $M$ be two lists of the same degree
and $N$ be any list. Then we have,
$$
\begin{array}{r c l}
\beta(L) > \beta(M) & \mbox{\rm iff} 
         & \beta(L \centerdot N) > \beta(M \centerdot N) \\
\beta(L) = \beta(M) & \mbox{\rm iff} 
         & \beta(L \centerdot N) = \beta(M \centerdot N).
\end{array}
$$
\end{lemma}

\begin{lemma} \label{l:switch_list}
Let $M$ and $N$ be any two lists.
Then
$\beta(M \centerdot N) = \beta(M^* \centerdot N) =
\beta(M \centerdot N^*)$
and
$\beta(M \centerdot N) = \beta(N \centerdot M)$.
\end{lemma}

\begin{lemma} \label{l:bullet_list}
Let 
$M = (M',m)$
and
$N = (n,N')$
be any two lists. Then
\[
(M) \centerdot (N) = (M',m-1,n,N') + (M',m,n-1,N') + 2 (M',m+n+1,N').
\]
\end{lemma}

\begin{lemma} \label{l:definition_S_list}
Let $M = (m_1,\ldots,m_i,\ldots,m_k)$ be any list. Then
\[
S(M) = \sum_{i=1}^k (m_1,\ldots,m_i-1,\ldots,m_k) +
        \sum_{i=1}^{k-1} (m_1,\ldots,m_i+m_{i+1}+1,\ldots,m_k).
\]
\end{lemma}

\begin{lemma} \label{l:alternate_S_list}
Let $M = (m_1,\ldots,m_i,\ldots,m_k)$ be any list. Then

$2 S(M) = (m_1-1,\ldots,m_i,\ldots,m_k) + 
                (m_1,\ldots,m_i,\ldots,m_k-1) +
                S'(M)$,
where

$S'(M) = \sum_{i=1}^{k-1} (m_1,\ldots,m_i) \centerdot (m_{i+1},\ldots,m_k)$.
\end{lemma}

\begin{lemma} \label{l:useful}
We make the following two useful observations.
\[
\begin{array}{r c l}
S(M,m,n,N) & = & (S(M,m),n,N) + (M,m+n+1,N) + (M,m,S(n,N)). \\
S'(M,m,n,N) & = & (S'(M,m),n,N) + (M,m) \centerdot (n,N) + (M,m,S'(n,N)).
\end{array}
\]
\end{lemma}

\begin{lemma} \label{l:identities_list}
Let $L,M$ and $N$ be any three lists. Then, we have,

$\beta(0,L,1,M) = \beta(0,L^*,1,M)$ 
and
$\beta(L,1,M,1,N) = \beta(L,1,M^*,1,N).$
\end{lemma}

\begin{lemma} \label{l:cc_d_list}
For any lists $M$ and $N$ and non-negative integers $m$ and $n$,
we have
$\beta(M,m,n,N) \geq \beta(M,m+n+2,N).$
\end{lemma}

\subsection{More applications} \label{ss:ma}

Though there are no new ideas in Section~\ref{s:ln},
the compactness of the list notation
allows us to make manipulations efficiently.
We illustrate this with some corollaries and exercises.
The reader, who is interested in later sections,
may skip ahead to Section~\ref{s:unimodal}.
Exercises 6 and 9 will be used later in 
Section~\ref{s:maxima}, where we tackle the problem
of locating the maxima. 
Other than that, 
no other results in this section will be used later.

\begin{corollary} \label{c:identity}
Let $M$ be of the form $(0,\ldots,0)$.
Also let $i,k$ be non-negative integers.
Then
\begin{equation} \label{e:identity}
\beta(i,M,i+k) - \beta(i,i+k,M) =
\beta(i+k-1,M,i+1) - \beta(i+k-1,i+1,M).
\end{equation}
\end{corollary}
\begin{proof}
We will show that equation~\eqref{e:identity} holds 
by induction on the number of zeroes in $M$.
When $M$ is empty, the result clearly holds.

By Lemma~\ref{l:switch_list} and the fact that $M=M^*$,
we have the two identities
\[
\begin{array}{c}
\beta [(i+1) \centerdot (M,i+k)] = \beta [(i+1) \centerdot (i+k,M)]\\
\beta [(i+k) \centerdot (M,i+1)] = \beta [(i+k) \centerdot (i+1,M)].
\end{array}
\]
Expanding the first identity using Lemma~\ref{l:bullet_list}
and regrouping terms,
we see that the left hand side of equation~\eqref{e:identity} is
\[
\beta(i+1,i+k-1,M) - 2 \beta(i+2,M',i+k) + 2 \beta(2i+k+2,M),
\]
where $M'$ has the same form as $M$ but with one zero less.
Similarly, expanding the second identity shows that 
the right hand side of equation~\eqref{e:identity} is
\[
\beta(i+k,i,M) - 2 \beta(i+k+1,M',i+1) + 2 \beta(2i+k+2,M).
\]
Now if we expand the identity
$\beta((i+k,i+1) \centerdot (M)) = \beta((i+1,i+k) \centerdot (M))$,
we reduce ourselves to proving equation~\eqref{e:identity},
but with $M$ replaced by $M'$.
Hence the result follows by induction.
\end{proof}

\noindent
As an example, 
for the Boolean lattice $B_9$,
Corollary~\ref{c:identity} gives us the following three identities.
\[
\begin{array}{l}
\beta(0,0,4) - \beta(0,4,0) = \beta(3,0,1) - \beta(3,1,0),\\
\beta(1,0,3) - \beta(1,3,0) = \beta(2,0,2) - \beta(2,2,0) \ \ \text{and}\\
\beta(2,0,0,0) - \beta(0,2,0,0) = \beta(1,0,0,1) - \beta(1,1,0,0).
\end{array}
\]

\begin{corollary} \label{c:ul}
Let $k$ be a non-negative integer and $K$ be any list. Then
\[
(k+1) \beta(k+2,K) \leq
\beta(0,k,K) \leq
(k+2) \beta(k+2,K).
\]
\end{corollary}
\begin{proof}
We show the second inequality.
The first inequality can be proved similarly.
By Lemma~\ref{l:reduction},
it is enough to show 
$2 \beta(S(0,k,K)) \leq
2(k+2) \beta(S(k+2,K))$.
Expand both sides using
Lemma~\ref{l:alternate_S_list}.
Hence we are reduced to showing that
\[
\beta((0) \centerdot (k,K)) + \beta((0,k) \centerdot (K))
\leq
(k+2)[\beta((k+2) \centerdot (K)) + \beta(k+1,K)].
\]
We know from Example~\ref{eg:compute} that 
$\beta(0,k) = k+1 = (k+1)\beta(k+2)$.
Applying Lemma~\ref{l:add_on_list}, we get
$\beta((0,k) \centerdot K) =
(k+1) \beta((k+2) \centerdot K)$.
Now Lemma~\ref{l:bullet_list} implies
\[
\beta((0,k) \centerdot K) \leq
(k+2) \beta((k+2) \centerdot K) -
\beta(k+1,K).
\]
Again by Lemma~\ref{l:bullet_list}, we have
$\beta((0) \centerdot (k,K)) =
\beta(0,k-1,K) + 2 \beta(k+1, K)$.
Using the induction hypothesis yields
\[
\beta((0) \centerdot (k,K)) \leq
(k+3) \beta(k+1,K).
\]
Adding the last two inequalities,
we obtain the desired result.
\end{proof}

We suggest some useful exercises, 
which the reader might want to try out.
Some of the inequalities that occur here are
very interesting.
They are all proved using induction,
the starting point being Lemma~\ref{l:reduction}.
We ask the reader to compare exercises 1 and 2,
exercises 3 and 4 
and also exercises 7 and 8.
They have the same flavour as Corollary~\ref{c:ul},
where we got upper and lower bounds 
for the $\beta$ value of a $\cd$-monomial.

\begin{exercise} \label{exercise_1}

$\beta(1,0,M) \geq \beta(0,1,M)$.

Hint: Use Lemma~\ref{l:alternate_S_list}. 
It is clear that the same proof also gives us
$\beta(k,0,M) \geq \beta(0,k,M)$.
For the most general result in this direction,
see Theorem~\ref{t:adjoining_balance}, part $(1)$.
\end{exercise}

\begin{exercise} \label{exercise_2}
$2 \beta(0,1,M) = \beta(1,0,M) + 2 \beta(3,M)$.

Hint: Since $\beta(0,0) = \beta(2)$,
Lemma~\ref{l:add_on_list} gives
$\beta((0,0) \centerdot (0,M)) = \beta((2) \centerdot (0,M))$.
Now use Lemma~\ref{l:bullet_list}.
\end{exercise}

\begin{exercise} \label{exercise_3}
$\beta(0,0,0,M) \geq \beta(2,0,M) + \beta(4,M)$.

Hint: Use Lemma~\ref{l:alternate_S_list}, 
along with the facts
$\beta(0,0,0) = \beta(2,0) + \beta(4)$ and
$2 \beta(0,1,M) \geq 3 \beta(3,M)$.
The second fact follows from the previous exercise
and Lemma~\ref{l:cc_d_list}.
\end{exercise}

\begin{exercise} \label{exercise_4}
$\beta(2,0,M) + 2 \beta(4,M) > \beta(0,0,0,M)$
for $M$ of the form $(0,\ldots,0)$.

Hint: We suggest that the reader
do the next two exercises, which imply this one.
\end{exercise}

\begin{exercise} \label{exercise_5}
$\beta(2,0,M) + 2 \beta(4,M) = \beta(1,1,M)$.

Hint: By Lemma~\ref{l:add_on_list}, we get
$\beta((3) \centerdot (0,M)) = 1/2 \beta((1,0) \centerdot (0,M)).$
Now use Lemma~\ref{l:bullet_list}.
\end{exercise}

\begin{exercise} \label{exercise_6}
$\beta(1,1,M) > \beta(0,0,0,M)$ 
for $M$ of the form $(0,\ldots,0)$.

Hint: Applying the usual method and using
$\beta(1,1) = \beta(0,0,0) + \beta(4)$,
this inequality reduces to the one in exercise~\ref{exercise_7}.
Since exercise~\ref{exercise_7} reduces to this exercise,
both statements can be proved by a joint induction.
\end{exercise}

\begin{exercise} \label{exercise_7}
$3 \beta(3,0,M) + 2 \beta(5,M) > \beta(1,0,0,M)$
for $M$ of the form $(0,\ldots,0)$.

Hint: Use the usual method 
along with the result of exercise~\ref{exercise_5}
to reduce to the previous exercise.
\end{exercise}

\begin{exercise} \label{exercise_8}
$\beta(1,0,0,M) \geq 3 \beta(3,0,M)$.

Hint: Use Lemma~\ref{l:alternate_S_list}, 
along with exercise~\ref{exercise_3} and the facts
$\beta((1,0) \centerdot (0,M)) = 2 \beta((3) \centerdot (0,M))$
and
$\beta(1,0,0) = 3 \beta(3,0)$.
\end{exercise}

\begin{exercise} \label{exercise_9}
$\beta(0,0,0,0,0,M) > \beta(0,1,1,0,M)$
for $M$ of the form $(0,\ldots,0)$.

Hint: Use Lemma~\ref{l:alternate_S_list} and the facts
$(0) \centerdot (1,1,0,M) = (0,1) \centerdot (1,0,M)$,
$\beta(0,1,1) = \beta(0,0,0,0) + \beta(6)$
and $\beta(1,0,0,M) \geq 2 \beta(0,3,M)$.
The first fact is true because
by Lemma~\ref{l:add_on_list},
both sides are equal to $2 (1,1,M',0,1)$,
where $M'$ has one zero less than $M$.
\end{exercise}

Restatements of exercises~\ref{exercise_6} and~\ref{exercise_9}
are given in equations~\eqref{e:ex6} and~\eqref{e:ex9} respectively
in Section~\ref{s:maxima}.

\section{Unimodal sequences} \label{s:unimodal}

In this section,
we study patterns of reverse unimodal sequences 
that arise in the $\cd$-index of the Boolean lattice.
For simplicity of notation, from now on, 
we will write
$(0^s)$ for the list $(\underbrace{0,\ldots,0}_{s})$.

\subsection{A unimodal sequence}
\label{ss:unimodal}

Recall the Euler numbers $E_n, n \geq 0$
defined by
$$\tan(x) + \sec(x) = \sum_{n \geq 0} E_n \cdot \frac{x^n}{n!}.$$
In other words, the odd and even Euler numbers 
are the tangent and secant numbers respectively.
In what follows, we will be dealing only with the tangent numbers.
We first recall an interesting inequality
involving the tangent numbers.
It is a special case of~\cite[Proposition 7.1]{\emahajan}.

\begin{proposition} \label{p:euler}
Let $a,b,c,d$ be non-negative odd integers, such that
$a + b = c + d = n$.
Then for $|a - b| > |c - d|$,
we have 
$$     \binom {n}{a} \cdot E_{a} \cdot E_{b}
            >
       \binom {n}{c} \cdot E_{c} \cdot E_{d}   . $$
\end{proposition}

Our interest in the odd Euler numbers comes from the fact
that they are related to the sequence
$\beta(\dv), \beta(\dv^2), \beta(\dv^3), \ldots$.
We prove the following lemma.

\begin{lemma} \label{l:deuler}
We have $2^n \beta(\dv^n) = E_{2n+1}.$
\end{lemma}
\begin{proof}
Consider the generating function
\[
P(x) = \sum_{n \geq 0} \frac{2^n \beta(\dv^n)}{(2n+1)!} \cdot x^{2n+1}.
\]
Note that $\sec(x)$ and $\tan(x)$ are even and odd functions 
respectively.
So by the definition of the Euler numbers,
the lemma is equivalent to showing that
$P(x) = \tan(x)$.
Observe that
\[
2 \beta(\dv^{n+1}) = 
\sum_{i=0}^n \beta(\dv^i \centerdot \dv^{n-i}) =
\sum_{i=0}^n \binom {2n+2}{2i+1} \beta(\dv^i) \beta(\dv^{n-i}).
\]
The first equality follows from Lemmas~\ref{l:reduction} and 
\ref{l:formula_useful}
and the second from Lemma~\ref{l:easy}. 
After elementary manipulations,
this gives us the differential equation
$P'(x) = 1 + P(x)^2$.
Using the initial conditions
$P(0)=0$ and $P'(0)=1$,
we conclude that $P(x)=\tan(x)$.
\end{proof}

\begin{remark}
Ehrenborg pointed out that the above lemma
is also a consequence of~\cite[Proposition 8.2]{\ber}.
\end{remark}

\begin{proposition} \label{p:unimodal}
The sequence 
$\beta(0^i,1,0^{n-i})$
for $0 \leq i \leq n$
is reverse unimodal in $i$.
\end{proposition}
In other words, as $i$ increases,
$\beta(0^i,1,0^{n-i})$ decreases till $i = [n/2]$
and then increases again.
Due to the symmetry property
$\beta(L) = \beta(L^*)$,
the proposition can be equivalently stated as follows.

Let $i,j,l,m$ be non-negative integers such that
$i+j = l+m = n$ and $|i-j| > |l-m|$. Then
$\beta(0^i,1,0^j) > \beta(0^l,1,0^m)$.

\begin{proof}
Note that by Lemma~\ref{l:bullet_list}, we obtain
$2 (0^i,1,0^j) = (0^{i+1}) \centerdot (0^{j+1}) =
\dv^i \centerdot \dv^j$ 
and by Lemma~\ref{l:easy}, we have
$\beta(\dv^i \centerdot \dv^j) = 
\binom {2i+2j+2}{2i+1} \beta(\dv^i) \beta(\dv^j).$
Hence equivalently,
we want to show that for $i,j,l,m$ as above
\[
\binom {2n+2}{2i+1} \beta(\dv^i) \beta(\dv^{j}) >
\binom {2n+2}{2l+1} \beta(\dv^l) \beta(\dv^{m}).
\]
This is a consequence of Proposition~\ref{p:euler} 
and Lemma~\ref{l:deuler}.
\end{proof}

\subsection{A general conjecture} \label{ss:gc}

Let $l$ be a positive integer
and $i$ and $j$ be non-negative integers such that $j > i$.
Let $L_{i,j}^k$ be the list 
of length $l$ given by
$(i,\ldots,i,j,i,\ldots,i)$,
where the letter $j$ appears 
in the $k$th position.
For simplicity of notation,
we have suppressed $l$ in our notation.
We are interested in the families of lists
where we fix $i,j$ and $l$ and
let $k$ vary from $1$ to $l$.
Define
$a_{i,j}^k = \beta(L_{i,j}^k)$.

\begin{lemma}
Let $k$ and $k'$ be positive integers. Then
for a fixed list length, we have
\[
\{a_{0,1}^k\} > \{a_{0,1}^{k'}\}
\ \text{iff} \ 
\{a_{0,2}^k\} > \{a_{0,2}^{k'}\}.
\]
\vanish{
$\{a_{0,1}^k\} > \{a_{0,1}^{k'}\}$
iff 
$\{a_{0,2}^k\} > \{a_{0,2}^{k'}\}$
}
\end{lemma}
\begin{proof}
The lemma follows from the following chain of equalities.
\[
\begin{array}{r c l l}
\beta((0) \centerdot (0^s,1,0^t)) & = & 
2 \beta(1,0^{s-1},1,0^t) & \\
& = & \beta((1,0^{s}) \centerdot (0^{t+1})) & \\
& = & \beta((0^{s},1) \centerdot (0^{t+1})) & \\
& = & \beta(0^{s+t+2}) + 2 \beta(0^s,2,0^t) & .\\
\end{array}
\]
The third equality follows from Lemma~\ref{l:switch_list}
while the remaining ones follow from Lemma~\ref{l:bullet_list}.
\end{proof}

The lemma, in particular, says that 
for a fixed list length $l$,
the sequence $\{a_{0,1}^k\}$ is reverse unimodal in $k$
iff
$\{a_{0,2}^k\}$ is reverse unimodal in $k$.
This gives us the following corollary
to Proposition~\ref{p:unimodal}.

\begin{corollary} \label{c:unimodal}
The sequence 
$\beta(0^i,2,0^{n-i})$
for $0 \leq i \leq n$
is reverse unimodal in $i$.
\end{corollary}

Motivated by the results so far,
we make a general conjecture.

\begin{conjecture} \label{conjecture_unimodal}
In the notation above,
for fixed $i,j$ and $l$,
the sequence $\{a_{i,j}^k\}$ is reverse unimodal in $k$.
\end{conjecture}

From Proposition~\ref{p:unimodal} and
Corollary~\ref{c:unimodal},
we know that the conjecture holds for the cases
$i=0,j=1$ and $i=0,j=2$ respectively.
The first step in the general conjecture,
namely,
$\{a_{i,j}^1\} > \{a_{i,j}^2\}$,
is a special case of
Theorem~\ref{t:adjoining_balance}, part $(1)$
in Section~\ref{s:balance}.
The remaining cases of the conjecture are open.

\begin{remark}
It might be possible to replace ``unimodal''
by ``log-concave''. Also there might be many other
families of such sequences that we have not accounted for.
There is a rich variety of methods for showing
that a sequence is log-concave or unimodal.
We refer the reader to the survey paper by Stanley~\cite{\stanleylog}.
\end{remark}

\section{Locating the maximum} \label{s:maxima}

We are now ready to answer the question that
was raised in item (3) of Section~\ref{ss:questions},
namely,
that of finding the $\cd$-monomial 
whose $\beta$ value is maximum.
We first review some of the facts that we need.
As in the previous section, we write
$(0^s)$ for the list $(\underbrace{0,\ldots,0}_{s})$.

A list with entries $0$ and $1$ can be written
(upto a power of $2$) as a product of lists
of the form $(0^s)$.
Also note that the list $(0)$ can appear
at most twice in the factorisation.
For example,
\[
(0,1,1,0,1) = 1/8 \ (0,0) \centerdot (0,0) \centerdot (0,0,0) \centerdot (0) = 
1/8 \ (0^2) \centerdot (0^2) \centerdot (0^3) \centerdot (0).
\]
This follows from the explicit description 
of the $\centerdot$ product given by Lemma~\ref{l:bullet_list}. 
Using this observation,
exercises~\ref{exercise_6} and \ref{exercise_9} 
can be restated as follows.
\begin{equation} \label{e:ex6}
\beta((0) \centerdot (0^2) \centerdot (0^s)) \geq 4 \ \beta(0^{s+2})
\quad \text{for} \quad s \geq 1.
\end{equation}
\begin{equation} \label{e:ex9}
4 \ \beta(0^{s+3}) > \beta((0^2) \centerdot (0^2) \centerdot (0^s))
\quad \text{for} \quad s \geq 1.
\end{equation}
And a special case of Proposition~\ref{p:unimodal} says that
\begin{equation} \label{e:eunimodal}
\beta((0) \centerdot (0^{n-1})) \geq \beta((0^i) \centerdot (0^{n-i}))
\quad \text{for} \quad i,n-i \geq 1.
\end{equation}
\begin{equation} \label{e:eeunimodal}
\beta((0^2) \centerdot (0^{n-2})) \geq \beta((0^i) \centerdot (0^{n-i}))
\quad \text{for} \quad i,n-i \geq 2.
\end{equation}
We now prove the main result of this section.

\begin{theorem} \label{t:maxima}
Among all {\bf cd}-monomials of a given degree, 
${\bf c} {\bf d}^{n} {\bf c}$ 
or
${\bf c} {\bf d} {\bf c} {\bf d}^{n} {\bf c}$
and
${\bf c} {\bf d}^{n} {\bf c} {\bf d} {\bf c}$
are the maxima, 
depending on whether the degree is even or odd. 
In the list notation, the maxima are 
$(1,0,\ldots,0,1)$ or 
$(1,1,0,\ldots,0,1)$ and
$(1,0,\ldots,0,1,1)$.
\end{theorem}
\begin{proof}
The idea of the proof is as follows.
Start with any list $M$.
Modify it to obtain a new list $M'$
such that $\beta(M') \geq \beta(M)$ 
holds. Now repeat the process on $M'$.
Continue this procedure till the modified list 
is one of the three lists in the theorem.

Note that the three lists in the theorem have factorisations 
$(0) \centerdot (0^s) \centerdot (0)$ and
$(0) \centerdot (0^2) \centerdot (0^s) \centerdot (0)$
and
$(0) \centerdot (0^s) \centerdot (0^2) \centerdot (0)$
respectively.
We know from Lemma~\ref{l:switch_list} that
the order of the factors does not change the $\beta$ value.
Now we enumerate our list modifications sequentially.
Modify the list $M$ so that

\begin{itemize}
\item
The entries in $M$ are either $0$ or $1$.
Hence $M$ has a factorisation into lists of the form $(0^s)$,
with $(0)$ appearing at most twice.

This is done by repeatedly applying 
Lemma~\ref{l:cc_d_list}.

\item
$M$ has entries $0$ and $1$ and it begins and ends with $1$.
In other words, the list $(0)$ appears exactly twice 
in the factorisation of $M$.

This is done by applying equation~\eqref{e:ex6} or \eqref{e:eunimodal},
whichever is appropriate.

\item
$M$ has the form $(1,\ldots,1,0,\ldots,0,1)$.
In other words, the factorisation of $M$ has the form
$(0) \centerdot (0^2) \centerdot \ldots \centerdot (0^2) \centerdot (0^s) \centerdot (0)$ for some $s \geq 2$.

This is done by repeatedly using equation~\eqref{e:eeunimodal}.

\item
$M$ is a maxima, i.e., 
the factorisation of $M$ has the form
$(0) \centerdot (0^s) \centerdot (0)$ or
$(0) \centerdot (0^2) \centerdot (0^s) \centerdot (0)$.

This is done by repeatedly using equation~\eqref{e:ex9}.
\end{itemize}
\end{proof}
\noindent
We illustrate the process described in the proof on two examples.
$$
\begin{array}{r c l l}
\beta(0^8) & < &
1/4 \ \beta((0) \centerdot (0^2) \centerdot (0^6)) & 
\quad \text{equation~\eqref{e:ex6}} \\
                       & < &
1/16 \ \beta((0) \centerdot (0^2) \centerdot (0^2) \centerdot (0^4) 
\centerdot (0)) & 
\quad \text{equation~\eqref{e:ex6}} \\
                       & < &
1/4 \ \beta((0) \centerdot (0^7) \centerdot (0)) & 
\quad \text{equation~\eqref{e:ex9}} \\
                       & = &
\beta(1,0^5,1). & 
\end{array}
$$
$$
\begin{array}{r c l l}
\beta(0^2,1,0,1,0^4,1) & = &
1/8 \ \beta((0^3) \centerdot (0^3) \centerdot (0^6) \centerdot (0)) & \\
                       & < &
1/8 \ \beta((0) \centerdot (0^5) \centerdot (0^6) \centerdot (0)) & 
\quad \text{equation~\eqref{e:eunimodal}} \\
                       & < &
1/8 \ \beta((0) \centerdot (0^2) \centerdot (0^9) \centerdot (0)) & 
\quad \text{equation~\eqref{e:eeunimodal}} \\
                       & = &
\beta(1,1,0^7,1). & 
\end{array}
$$

\section{The balance inequalities}  \label{s:balance}

In this section,
we study inequalities that involve balancing of $\cd$-monomials.
The motivation for these considerations comes from
similar inequalities for the $\ab$-monomials
that were proved in~\cite{\emahajan}.
The intuitive connection between the two situations
is given in Appendix~\ref{s:connection}.

\subsection{The balancing of a {\bf cd}-monomial}

Let $m_1,m_2,n_1,n_2$ be non-negative integers such that
$m_1 + n_1 = m_2 + n_2$.
We say that a pair
$(m_1,n_1)$ is \emph{better balanced} than a pair 
$(m_2,n_2)$ if 
$|m_1 - n_1| \leq |m_2 - n_2|$.
And we say that it is \emph{strictly better balanced}
if the inequality is strict.

We check that under this condition, 
we can pair off the terms
$(m_1-1,n_1)$ and $(m_1,n_1-1)$ with the terms
$(m_2-1,n_2)$ and $(m_2,n_2-1)$,
not necessarily in the same order,
such that 
the same condition still holds for each of the two pairs.
We refer to this as the {\it reduction} property. 

\begin{lemma} \label{l:balance_step}
Let $(m_1,n_1)$ be a pair that is better balanced than the pair 
$(m_2,n_2)$. Then 
$\beta(m_1,n_1) \geq \beta(m_2,n_2)$ with equality iff
$|m_1 - n_1| = |m_2 - n_2|$.
\end{lemma}
\begin{proof}
We prove the result by induction on
$m_1 + n_1 = m_2 + n_2$.
By the reduction property and the induction hypothesis, 
we obtain
\[
\beta(m_1-1,n_1) + \beta(m_1,n_1-1) 
\geq \beta(m_2-1,n_2) + \beta(m_2,n_2-1).
\]
Adding $\beta(m_1+n_1+1) = \beta(m_2+n_2+1)$ to both sides,
we get
$\beta(S(m_1,n_1)) \geq \beta(S(m_2,n_2))$, which by
Lemma~\ref{l:reduction},
yields the desired result.
\end{proof}

\begin{remark}
This lemma also follows from the formula
$\beta(m,n) = {\binom {m+n+2}{m+1}} - 1$,
which we wrote in Example~\ref{eg:compute}. 
However, we prefer the non-computational proof
above since it illustrates our basic technique.
The key idea is that applying $S$ does not change the $\beta$ value
and $S$ has degree $-1$. 
In what follows, we will also use the description of $S$
given by Lemma~\ref{l:alternate_S_list}.
It involves the map $S'$,
which is also of degree $-1$.
\end{remark}

Now we state the main result of this section.

\begin{theorem} \label{t:adjoining_balance}
Let $m$,$n$ be non-negative 
integers such that
$n > m$.
Let $(m_1,n_1)$ be a pair that is strictly better balanced than the pair 
$(m_2,n_2)$. 
Also let 
$M$ and $N$ be any two lists. 
Then we have

\begin{itemize}
\item[(1)]
$A(r,l) : \: \: \beta(M,m,n) \geq \beta(M,n,m),$ with equality if $M$ is empty.

\item[(2)]
$B(r,l) : \: \: \beta(M,m_1,n_1,N) > \beta(M,m_2,n_2,N)$. 
\end{itemize}

\noindent
{\rm
The letters $r$ and $l$ denote 
the degree and length of the lists
that appear in the two statements.}
\end{theorem}
\begin{proof}
We prove parts $(1)$ and $(2)$ of the theorem using a joint
induction. 
The induction is on $r$ and $l$ and is divided in three steps.
The first step is the induction basis.
The next two are the induction steps
for parts (1) and (2) respectively.

\subsection*{(i)}
The induction basis for part (1) is
the statement $A(r,2)$, which just says 
$\beta(m,n) = \beta(n,m)$.
For part (2), it is the statement $B(r,2)$,
which says 
$\beta(m_1,n_1) > \beta(m_2,n_2)$.
This is true by Lemma~\ref{l:balance_step}.

\subsection*{(ii)}
$A(< r, \leq l)$ and $B(< r, < l)$ implies $A(r,l)$.
%To show $\beta(M,m,n) \geq \beta(M,n,m)$.

Set $M = (K,k)$. 
Then Lemma~\ref{l:definition_S_list} and
the statement $A(< r, \leq l)$  gives
\[
\beta(K,k,S(m,n)) \geq \beta(K,k,S(n,m))
\ \ \text{and} \ \ 
\beta(S(K,k),m,n) \geq \beta(S(K,k),n,m).
\]
Also statement $B(< r, < l)$ gives
$\beta(K,k+m+1,n)) \geq \beta(K,k+n+1,m))$.
Summing up the last three inequalities
and using Lemma~\ref{l:useful}, we obtain
\[
\beta(S(K,k,m,n)) \geq \beta(S(K,k,n,m)).
\]
Now by Lemma~\ref{l:reduction},
we get
$\beta(M,m,n) \geq \beta(M,n,m)$,
which is the statement $A(r,l)$.

\subsection*{(iii)}
$A(r,l)$ and $B(< r, \leq l)$ and $B(r, < l)$ implies $B(r,l)$.
%To show $\beta(M,m_1,n_1,N) > \beta(M,m_2,n_2,N)$. 

We split this step into two cases.

\medskip
\noindent
{\it Case 1:}
$M$ and $N$ are both non-empty.

Using Lemma~\ref{l:add_on_list},
the definition of the map $S'$
given by Lemma~\ref{l:alternate_S_list}
and
the statement $B(< r, < l)$, we get
\[
\begin{array}{c}
\beta(S'(M,m_1),n_1,N) > \beta(S'(M,m_2),n_2,N). \\
\beta(M,m_1,S'(n_1,N)) > \beta(M,m_2,S'(n_2,N)).
\end{array}
\] 
Also by the reduction property,
Lemma~\ref{l:bullet_list} and 
the statement $B(< r,l)$, we get
$\beta((M,m_1) \centerdot (n_1,N)) > \beta((M,m_2) \centerdot (n_2,N))$. 
Summing up the last three inequalities
and using Lemma~\ref{l:useful}, we obtain
\[
\beta(S'(M,m_1,n_1,N)) > \beta(S'(M,m_2,n_2,N)).
\]
Since $M$ and $N$ are non-empty, by $B(<r,l)$ we also have
\[
\begin{array}{c}
\beta(M-1,m_1,n_1,N) > \beta(M-1,m_2,n_2,N), \\
\beta(M,m_1,n_1,N-1) > \beta(M,m_2,n_2,N-1),
\end{array}
\]
where $(M-1)$ denotes one deleted from the first entry of $M$
and $(N-1)$ denotes one deleted from the last entry of $N$.
Adding the last three inequalities, and using
Lemma~\ref{l:alternate_S_list},
we obtain
$\beta(S(M,m_1,n_1,N)) > \beta(S(M,m_2,n_2,N))$.
Applying Lemma~\ref{l:reduction} gives statement $B(r,l)$.

\medskip
\noindent
{\it Case 2:}
Either $M$ or $N$ is empty.

Due to the symmetry property $\beta(L) = \beta(L^*)$,
we may assume that $M$ is non-empty and $N$ is empty. 
Now repeat the above argument. 
The only step that requires care is the inequality
that involves $N-1$.
Since $N$ is empty,
we are required to prove
$\beta(M,m_1,n_1-1) > \beta(M,m_2,n_2-1)$.
In most cases, 
the pair $(m_1,n_1-1)$ is better balanced 
than the pair $(m_2,n_2-1)$.
And hence applying $B(<r,l)$ completes the proof, as before.

The only case when it
fails to work is when
$m_1 \geq n_1$ and
$m_2 = n_1 - 1$ and $n_2 = m_1 + 1$.
In other words, we want to show the following
special case of $B(r,l)$.
\[
\beta(M,m,n) > \beta(M,n-1,m+1) 
\ \ \text{for} \ \ 
m \geq n.
\]
Set $M = (K,k)$. 
By 
Lemma~\ref{l:switch_list}, we have
$\beta((K,k) \centerdot (m+1,n)) = \beta((K,k) \centerdot (n,m+1))$.
Expand both sides using 
Lemma~\ref{l:bullet_list}. 
This gives us three terms on either side.
Using statements $A(r,l)$ and $B(r,< l)$ respectively,
two of the three terms 
can be compared as follows.
\[
\begin{array}{c}
\beta(K,k-1,m+1,n) < \beta(K,k-1,n,m+1). \\
\beta(K,k+m+2,n) < \beta(K,k+n+1,m+1).
\end{array}
\]
Hence for the remaining term, we obtain the inequality
$\beta(K,k,m,n) > \beta(K,k,n-1,m+1)$,
which is what we wanted to show.

\end{proof}
\noindent
Motivated by the previous theorem, 
we make the following conjectures.

\begin{conjecture} \label{c:skewed}
Let $m$,$n$ be non-negative 
integers such that
$n > m$.
Also let 
$L$ and $M$ be any two lists. Then
$\beta(M,m,L,n) \geq \beta(M,n,L,m),$ if $M$ is non-empty.
\end{conjecture}

\begin{conjecture} \label{c:balance_step}
Let $(m_1,n_1)$ be a pair that is strictly better balanced than the pair 
$(m_2,n_2)$. 
And let 
$L,M$ and $N$ be any three lists. 
Then we get
$\beta(M,m_1,L,n_1,N) > \beta(M,m_2,L,n_2,N)$. 
\end{conjecture}
\noindent
To state the next conjecture,
we require the notion of a balanced list.
A list $B$ is called {\it balanced} if its entries
are either $k$ or $k+1$ for some non-negative integer $k$.

\begin{conjecture} \label{c:balance}
Let $L$, $M$ and $L^{\prime}$ be three lists.
Then there exists a balanced list
$B$ of the same degree and length as $M$
such that $\beta(L,B,L^{\prime}) \geq \beta(L,M,L^{\prime})$.
\end{conjecture}

When $L$ is empty, conjectures~\ref{c:skewed} and 
\ref{c:balance_step} reduce to 
Theorem~\ref{t:adjoining_balance}.
Also conjecture~\ref{c:balance_step} implies
conjecture~\ref{c:balance}.

\subsection{A sufficient condition}

We are mainly interested in conjecture~\ref{c:balance_step}.
For the remainder of this section, 
we prove some of its special cases,
which are not accounted for by
Theorem~\ref{t:adjoining_balance}.
At the end of the section, we also give a sufficient
condition for its validity; see Theorem~\ref{t:suff_cond}.

\begin{lemma}
Let $(m_1,n_1)$ be a pair that is strictly better balanced than the pair 
$(m_2,n_2)$. 
Also let $M$ and $N$ be any two lists. Then
$\beta((m_1,M) \centerdot (N,n_1)) >
\beta((m_2,M) \centerdot (N,n_2))$.
\end{lemma}
\begin{proof}
The proof follows from the following chain of comparisons.

$\beta((m_1,M) \centerdot (N,n_1)) =
\beta((M^*,m_1) \centerdot (n_1,N^*)) >
\beta((M^*,m_2) \centerdot (n_2,N^*)) =
\beta((m_1,M) \centerdot (N,n_1))$.

The first and third equality follows from 
Lemma~\ref{l:switch_list}.
For the second inequality,
we expand both sides using Lemma~\ref{l:bullet_list}
and then use the reduction property
and Theorem~\ref{t:adjoining_balance}.
\end{proof}

\begin{lemma} \label{l:end_balance}
Let $(m_1,n_1)$ be a pair that is strictly better balanced than the pair 
$(m_2,n_2)$. 
Also let $L$ be any list. Then
$\beta(m_1,L,n_1) >
\beta(m_2,L,n_2)$.
\end{lemma}
\begin{proof}
We do an induction on the degree of the lists.
By the previous lemma and the definition of the map $S'$
given by Lemma~\ref{l:alternate_S_list},
we have
$\beta(S'(m_1,L,n_1)) >
\beta(S'(m_2,L,n_2))$.
Also by the reduction property and induction,
we get
\[
\beta(m_1-1,L,n_1) + \beta(m_1,L,n_1-1) 
\geq \beta(m_2-1,L,n_2) + \beta(m_2,L,n_2-1).
\]
Adding up the two inequalities and 
again using Lemma~\ref{l:alternate_S_list},
we get the inequality
\[
\beta(S(m_1,L,n_1)) >
\beta(S(m_2,L,n_2)).
\]
The result now follows from Lemma~\ref{l:reduction}.
\end{proof}

Repeated use of
Theorem~\ref{t:adjoining_balance}, part $(2)$ 
and
Lemma~\ref{l:end_balance}
proves the following.

\begin{corollary}
Conjecture~\ref{c:balance}
is correct in the special
case when
$L,L'$ are empty and $M$ is a list whose
length is at most three.
\end{corollary}

Next we prove two results that have the same flavour
as the previous two lemmas.

\begin{lemma}
Let $(m_1,n_1)$ be a pair that is strictly better balanced than the pair 
$(m_2,n_2)$. 
Let $n_2 < m_2$.
Also let 
$M$ and $L$ be any lists. Then we have
\[
\beta((M,m_1,L) \centerdot (n_1)) >
\beta((M,m_2,L) \centerdot (n_2)).
\]
\end{lemma}
\begin{proof}
We prove the result by induction.

\medskip
\noindent
{\it Induction basis:}
Either $M$ or $L$ is empty.

By Lemma~\ref{l:switch_list},
we may assume that $L$ is empty.
By Theorem~\ref{t:adjoining_balance} and the reduction property,
\[
\beta(M,m_1-1,n_1) + \beta(M,m_1,n_1-1) 
\geq \beta(M,m_2-1,n_2) + \beta(M,m_2,n_2-1).
\]
This inequality and Lemma~\ref{l:bullet_list} imply that
\[
\beta((M,m_1) \centerdot (n_1)) >
\beta((M,m_2) \centerdot (n_2)),
\]
which is what we wanted to show.

\medskip
\noindent
{\it Induction step:}
$M$ and $L$ are both non-empty.

By induction and Lemmas~\ref{l:switch_list},
\ref{l:bullet_list} and \ref{l:alternate_S_list},
we obtain
\[
\beta(S(M,m_1,L) \centerdot (n_1)) >
\beta(S(M,m_2,L) \centerdot (n_2)).
\]
Since $n_2 < m_2$,
the pair $(m_1,n_1-1)$ is strictly better balanced than 
the pair $(m_2,n_2-1)$.
Hence we apply induction to get
\[
\beta((M,m_1,L) \centerdot (n_1-1)) >
\beta((M,m_2,L) \centerdot (n_2-1)).
\]
Note that $S(n) = (n-1)$ for any $n \geq 1$.
Now adding the last two inequalities and applying Theorem~\ref{t:derivation},
we get
\[
\beta(S((M,m_1,L) \centerdot (n_1))) >
\beta(S((M,m_2,L) \centerdot (n_2))).
\]
Hence the result follows from Lemma~\ref{l:reduction}.
\end{proof}

\begin{lemma}
Let $(m_1,n_1)$ be a pair that is strictly better balanced than the pair 
$(m_2,n_2)$. 
Let $n_2 < m_2$.
Also let 
$L,M$ be any lists. Then
$\beta(M,m_1,L,n_1) > \beta(M,m_2,L,n_2)$.
\end{lemma}
\begin{proof}
The proof proceeds by induction.
If $M$ is empty, then the result holds by
Lemma~\ref{l:end_balance}.
This is the induction basis.

For the induction step,
we assume that $M$ is non-empty.
By induction,
the definition of the map $S'$
given by Lemma~\ref{l:alternate_S_list}
and Lemmas~\ref{l:switch_list} and~\ref{l:bullet_list}, 
we have
\[
\beta(S'(M,m_1,L),n_1) > \beta(S'(M,m_2,L),n_2).
\]
And by the previous lemma, we have
$\beta((M,m_1,L) \centerdot (n_1)) > \beta((M,m_2,L) \centerdot (n_2))$.
Adding the last two inequalities, we obtain
\[
\beta(S'(M,m_1,L,n_1)) > \beta(S'(M,m_2,L,n_2)).
\]
In addition, by induction, we also have
\[
\begin{array}{c}
\beta(M-1,m_1,L,n_1) > \beta(M-1,m_2,L,n_2).\\
\beta(M,m_1,L,n_1-1) > \beta(M,m_2,L,n_2-1).
\end{array}
\]
For the second inequality,
since $n_2 < m_2$,
the pair $(m_1,n_1-1)$ is strictly better balanced than 
the pair $(m_2,n_2-1)$.
Adding the last three inequalities
and using Lemmas~\ref{l:alternate_S_list} and \ref{l:reduction},
we get the conclusion of the lemma.
\end{proof}

\begin{remark}
In the previous two lemmas,
we may replace the condition $n_2 < m_2$
by the weaker condition that
$(m_1,n_1-1)$ is strictly better balanced than 
$(m_2,n_2-1)$.
\end{remark}
Encouraged by our success, 
let us try to prove Conjecture~\ref{c:balance_step} by induction.
We may assume that $L$ is non-empty.
Now we have three cases.

\medskip
\noindent
{\it Case 1:}
$M$ and $N$ are both empty.

This follows directly from Lemma~\ref{l:end_balance}.

\medskip
\noindent
{\it Case 2:}
$M$ and $N$ are both non-empty.

This case is again easy. We imitate the proof of 
Lemma~\ref{l:end_balance}. 
Note that we are relying on induction.

\medskip
\noindent
{\it Case 3:}
$M$ is non-empty and $N$ is empty.

In view of the previous lemma,
if we assume that 
$(m_1,n_1-1)$ is strictly better balanced than $(m_2,n_2-1)$,
then we have no trouble.
The only case for which this assumption does not work is
when
$m_1 \geq n_1$ and 
$m_2 = n_1 - 1$ and $n_2 = m_1 + 1$.

Therefore, we have the following sufficient condition for 
Conjecture 2 to hold.

\begin{theorem} \label{t:suff_cond}
Conjecture~\ref{c:balance_step} is true 
if it holds in the following special case.

Let $m$ and $n$ be non-negative integers such that
$m \geq n$. 
Also let $L$ and $M$ be any lists. 
Then
$\beta(M,m,L,n) > \beta(M,n-1,L,m+1)$.
\end{theorem}

The same problem as above arose 
while proving Theorem~\ref{t:adjoining_balance},
but we managed to deal with it there;
see Case $2$ in step (iii) of its proof.

\vanish{
\begin{theorem}
Let $(m_1,n_1)$ be a pair that is strictly better balanced than the pair 
$(m_2,n_2)$. 
Let $n_2 < m_2$.
Also let 
$K,L,M$ be any lists. Then

\begin{itemize}
\item[(1)] 
$\beta(M,m_1,L,n_1) > \beta(M,m_2,L,n_2)$, 

\item[(2)]  
$\beta(M,m_1,K,k,n_1) + \beta(M,m_1,K,n_1+k+2) > 
\beta(M,m_2,K,k,n_2) + \\
\beta(M,m_2,K,n_2+k+2)$.
\end{itemize}

Here $n_2$ is allowed to be $-1$ in statement (2).
\end{theorem}
\begin{proof}
The proof uses the same method as that for
Theorem~\ref{t:adjoining_balance}.
We give a sketch.

By Theorem~\ref{t:adjoining_balance},
statement (1) is true when the list $L$ is empty.
This provides the induction basis.
Denote statements (1) and (2) by $A(r)$ and $B(r)$ respectively.
We complete the induction step in 2 parts.

\begin{itemize}
\item[(i)] 
$A(<r), B(<r)$ implies $A(r)$.

Put $L = (K,k)$.
If $k > 0$ then $B(<r)$ implies
$\beta(M,m_1,K,k-1,n_1) + \beta(M,m_1,K,n_1+k+1) > 
\beta(M,m_2,K,k-1,n_2) + \beta(M,m_2,K,n_2+k+1)$.
This statement along with $A(<r)$ yields
$\beta(S(M,m_1,L,n_1)) > \beta(S(M,m_2,L,n_2))$. 

If $k = 0$ then put $K = (K',k')$.
Now use the previous lemma to obtain
$\beta(M,m_1,K',k'+1,n_1) + \beta(M,m_1,K',k',n_1+1) > 
\beta(M,m_2,K',k'-1,n_2) + \beta(M,m_2,K',k',n_2+1)$.
This statement along with $A(<r)$ yields the desired result.

\item[(ii)] 
$A(<r), B(<r)$ implies $B(r)$.

We essentially repeat the same argument as in item (i).
\end{itemize}
\end{proof}
}

\section{Concluding remarks} \label{s:conc_rmks}

We conclude with some comments and problems for further study.

\subsection{Divisibility properties}

As was mentioned in item (5) of Section~\ref{ss:questions},
many of the $\beta$ values are divisible by $1001$.
This phenomenon first occurs for 
$B_{13}$.
Upto list reversal and the identities provided by
Lemma~\ref{l:identities_list},
we provide a complete list of all
{\bf cd}-monomials of degree 12 
whose coefficients are divisible by 1001.
$$
\begin{array}{c c c}

\beta(6,1,1) = 5005 &

\beta(1,1,2,2) = 140140 &

\beta(2,1,1,2) = 162162 \\

\beta(3,1,1,1) = 120120 &

\beta(1,1,3,1) = 90090 &

\beta(2,1,3,0) = 54054 \\

\beta(1,1,0,4) = 50050 &

\beta(0,0,1,3,0) = 72072 &

\beta(1,1,1,1,0) = 300300 \\

\beta(2,0,0,1,1) = 260260 &

\beta(1,1,1,0,1) = 360360 &

\beta(2,1,0,1,0) = 216216 \\

\beta(0,1,0,1,0,0) = 288288 &&
\end{array}
$$

\vanish{
For $B_{14}$, 

$$
\begin{array}{c c c}
\beta(7,1,1) = 8008 &

\beta(6,1,2) = 21021 &

\beta(5,3,1) = 29029 \\

\beta(5,1,3) = 36036 &

\beta(3,5,1) = 23023 &

\beta(3,3,3) = 72072 \\

\beta(2,1,6) = 21021 &

\beta(1,7,1) = 4004 &

\beta(6,1,0,0) = 28028 \\

\beta(5,1,1,0) = 90090 &

\beta(5,1,0,1) = 108108 & 
\end {array}
$$
}
\noindent
This phenomenon continues for $B_{14}$,
where there are many more {\bf cd}-monomials
with this property.
We did not look at any data beyond rank $14$,
but we expect this behaviour to continue
and hence in need of some explanation.

\subsection{Recursions for the Boolean lattice} \label{ss:recursion}

Purtill~\cite{\purtill} gave the first recursion 
that showed that the Boolean lattice
had a {\bf cd}-index with positive coefficients.
$$
\Psi(B_{n+1}) = \cv \Psi(B_n) + \sum_{i=1}^{n-1} 
\binom {n-1}{i} \Psi(B_i) \dv \Psi(B_{n-i}).
$$
This recursion has a dual; the sum of the two gives a 
more symmetric recursion.
We did not make any use of these recursions in this paper. 
Instead, we worked with a certain derivation.

The Boolean lattice has a q-analogue,
namely the lattice of subspaces of a $n$ dimensional vector space
over the finite field $F_q$.
This lattice is usually denoted by $L_n$.
Then the {\bf ab}-index of the lattice of subspaces
satisfies the following recursion.
$$
\Psi(L_{n+1}) = (\av + q^n \bv) \Psi(L_n) + \sum_{i=1}^{n-1} 
\binom {n-1}{i}_q \Psi(L_i) (q^n \av \bv + q^i \bv \av) \Psi(B_{n-i}),
$$
with
$\Psi(L_1) = 1$,
$\Psi(L_2) = \av + q \bv$ and so on.
From this recursion, it looks unlikely that there is a nice
q-version of the {\bf cd}-index.
We also note that the expression that we have written down
in not unique.
For instance, this recursion also has a dual version; 
the sum of the two recursions then gives a third one.
These three recursions give three distinct ways of expressing 
the {\bf ab}-index of $\Psi(L_{n+1})$.

\subsection{An algebraic perspective}

Jointly with Marcelo Aguiar, a part of this paper has now been put in a
more algebraic context. 
The algebraic approach shows that
the existence of the coderivation $\hat G$ on $k \langle \av, \bv \rangle$
can also be derived from a certain universal property 
of the coalgebra $k \langle \av, \bv \rangle$.
It also gives an algebraic proof of the
recursions involving the {\bf ab}-index of the Boolean lattice 
and the lattice of subspaces
written in Section~\ref{ss:recursion}. 

\appendix
\section{Connection between the {\bf ab} and the {\bf cd}-index}
\label{s:connection}

In this section,
we point out some analogies between the results 
obtained in this paper and those in~\cite{\emahajan}
and give an intuitive explanation of why they occur.

We define a map
$\omega : k \langle \cv,\dv \rangle \rightarrow k \langle \av,\bv \rangle$.
For any {\bf cd}-monomial $v$,
let $\omega(v)$ be as follows.

Replace
every odd occurrence of $\dv$ in $v$ by $\av \bv$
and every even occurrence of $\dv$ in $v$ by $\bv \av$.
If the first $\dv$ to the right of a given $\cv$ in $v$ 
has an odd occurrence, 
then replace that $\cv$ by a $\av$,
else replace it by a $\bv$.
For example,
\[
\omega(\cv \dv \cv) = \av \av \bv \bv, \quad
\omega(\cv \dv \dv) = \av \av \bv \bv \av.
\]
The map $\omega$ is one-to-one and its image consists of those 
{\bf ab}-monomials, 
which begin with an $\av$
and which do not contain either
$\av \bv \av$ or $\bv \av \bv$ as a substring.
We will call such {\bf ab}-monomials {\it valid}.
%Extend $\omega$ to $k \langle \cv, \dv \rangle$ by linearity.

We now define a partial order on the set of all
{\bf cd}-monomials of a given degree as follows.

$v$ covers $u$ if 
$v$ may be obtained from $u$
by replacing an occurrence of $\cv^2$ in $u$ by a $\dv$.

The poset so defined is graded,
the rank of an element $v$ being the
number of occurrences of $\dv$ in $v$.
We denote the rank function by $\rho$.

We may transfer this partial order to the set of all
valid {\bf ab}-monomials of the same degree,
since the two sets are in bijection with each other.
This partial order may be described as follows.

$z$ covers $y$ if for some {\bf ab}-monomials 
$y_1$ and $y_2$, which end and begin respectively
with the same letter, we have $y = y_1 y_2$ 
and $z = y_1 {\overline{y}_2}$.
Here
${\overline{y}}_2$ is the {\bf ab}-monomial obtained from
$y_2$ by replacing an $\av$ by a $\bv$ and vice-versa.

Note that this partial order makes sense
for all $\ab$-monomials, not just for the valid ones.

For $v$, a {\bf cd}-monomial of degree $n$,
we have defined
$\beta(v)$ to be the coefficient of 
$v$ in $\Psi_{B_{n+1}}(\cv,\dv)$.
Similarly, for $y$, an {\bf ab}-monomial of degree $n$,
we define
$\beta(y)$ to be the coefficient of 
$y$ in $\Psi_{B_{n+1}}(\av,\bv)$.

\begin{lemma}
Let $v$ be any {\bf cd}-monomial. Then
\[
\beta(v) = \sum_{u \leq v} (-1)^{\rho(v)-\rho(u)} \beta(\omega(u)).
\]
\end{lemma}

The lemma follows directly from the definition of the
{\bf cd}-index and so we omit the proof.

\begin{remark}
We have stated this lemma only for the Boolean lattice.
But it holds for any poset that has a {\bf cd}-index.
It also shows that for an Eulerian poset
many entries of the flag $h$-vector (or the flag $f$-vector)
are redundant. The resulting linear relations are
the so called Dehn-Sommerville relations; see 
Theorem~\ref{t:Dehn_Sommerville}.
\end{remark}

Now we recall~\cite[Lemma 3.9]{\emahajan}
which says that
for any two {\bf ab}-monomials $y$ and $z$,
 inequality
$z \geq y$ implies $\beta(z) \geq \beta(y)$.
This has been referred to as the alternating property
in~\cite{\ereaddyr}.
It implies that in the alternating sum 
that occurs in the above lemma, 
the term with the largest magnitude is $\beta(\omega(v))$.
This gives us some reason to believe that 
if we linearly order the {\bf ab}-monomials
and linearly order the {\bf cd}-monomials by their
$\beta$ values, then the map $\omega$ would
respect this order to a large extent.
This is the intuition that led us to expect
$\cd$-analogues.
We point out three analogies.

\subsection*{(1)}
We have already noted that 
for any two {\bf ab}-monomials $y$ and $z$,
$z \geq y$ implies $\beta(z) \geq \beta(y)$.
Using the map $\omega$,
we expect the inequality
$\beta(u \dv v) \geq \beta(u \cv^2 v)$.
This result was obtained in Lemma~\ref{l:cc_d}.

\subsection*{(2)}
An {\bf ab}-monomial (that begins with an $\av$)
can be written uniquely written in the form
$\av^{m_1} \bv^{m_2} \av^{m_3} \ldots$,
where 
$m_1,m_2,\ldots,m_k$
are positive integers.
Hence we may represent it by the list
$(m_1,m_2,\ldots,m_k).$
This is the list notation for {\bf ab}-monomials
that was used in~\cite{\emahajan}.
This does not quite coincide under the map 
$\omega$ with our list notation for {\bf cd}-monomials,
but it is quite close.
Namely 
$\omega(m_1,m_2,\ldots,m_k) = (m_1+1,m_2+2,\ldots,m_{k-1}+2,m_k+1)$.

Next we recall the balance inequalities for {\bf ab}-monomials
that were shown in~\cite{\emahajan}.

\medskip
\begin{proposition}~\cite[Corollary 6.5]{\emahajan}
Let $(m_1,n_1)$ be a pair that is strictly better balanced than the pair 
$(m_2,n_2)$. 
Let $P$ be a palindrome and let $M$ and $N$ be any two lists.
Then
$\beta(M,m_1,P,n_1,N) > \beta(M,m_2,P,n_2,N)$. 
\end{proposition}

\medskip
\begin{proposition}~\cite[Theorem 6.7]{\emahajan}
Let $L$, $M$ and $L^{\prime}$ be three lists.
Then there exists a balanced list
$B$ of the same degree and length as $M$
such that $\beta(L,B,L^{\prime}) \geq \beta(L,M,L^{\prime})$.
\end{proposition}

The above result was originally conjectured by Gessel.
We proved the {\bf cd}-analogue of the first result for
the special case when $P$ is the empty list
(see Theorem~\ref{t:adjoining_balance})
and have conjectured the analogue for the second result
(see Conjecture~\ref{c:balance}).

\subsection*{(3)}
We recall some notation from Section~\ref{ss:gc}.
Let $L_{i,j}^k$ be the list 
of length $l$ given by
$(i,\ldots,i,j,i,\ldots,i)$,
where the letter $j$ appears 
in the $k$th position.
We are interested in the families of lists
where we fix $i,j$ and $l$ and
let $k$ vary from $1$ to $l$.
Define
$a_{i,j}^k = \beta(L_{i,j}^k)$.

If we think of the lists as $\ab$-monomials
then the numbers 
$a_{i,j}^k$ for $1 \leq k \leq l$ 
display a very intricate pattern as shown below.
\[
\begin{minipage}{2 in}
\vspace{.1in}
\begin{center}
\psfrag{a}{\Huge $a_{i,j}^k$}
\psfrag{k}{\Huge $k$}
\psfrag{1}{\Huge $1$}
\psfrag{2}{\Huge $2$}
\psfrag{3}{\Huge $3$}
\psfrag{4}{\Huge $4$}
\psfrag{5}{\Huge $5$}
\psfrag{6}{\Huge $6$}
\psfrag{7}{\Huge $7$}
\psfrag{8}{\Huge $8$}
\resizebox{2 in}{1 in}{\includegraphics{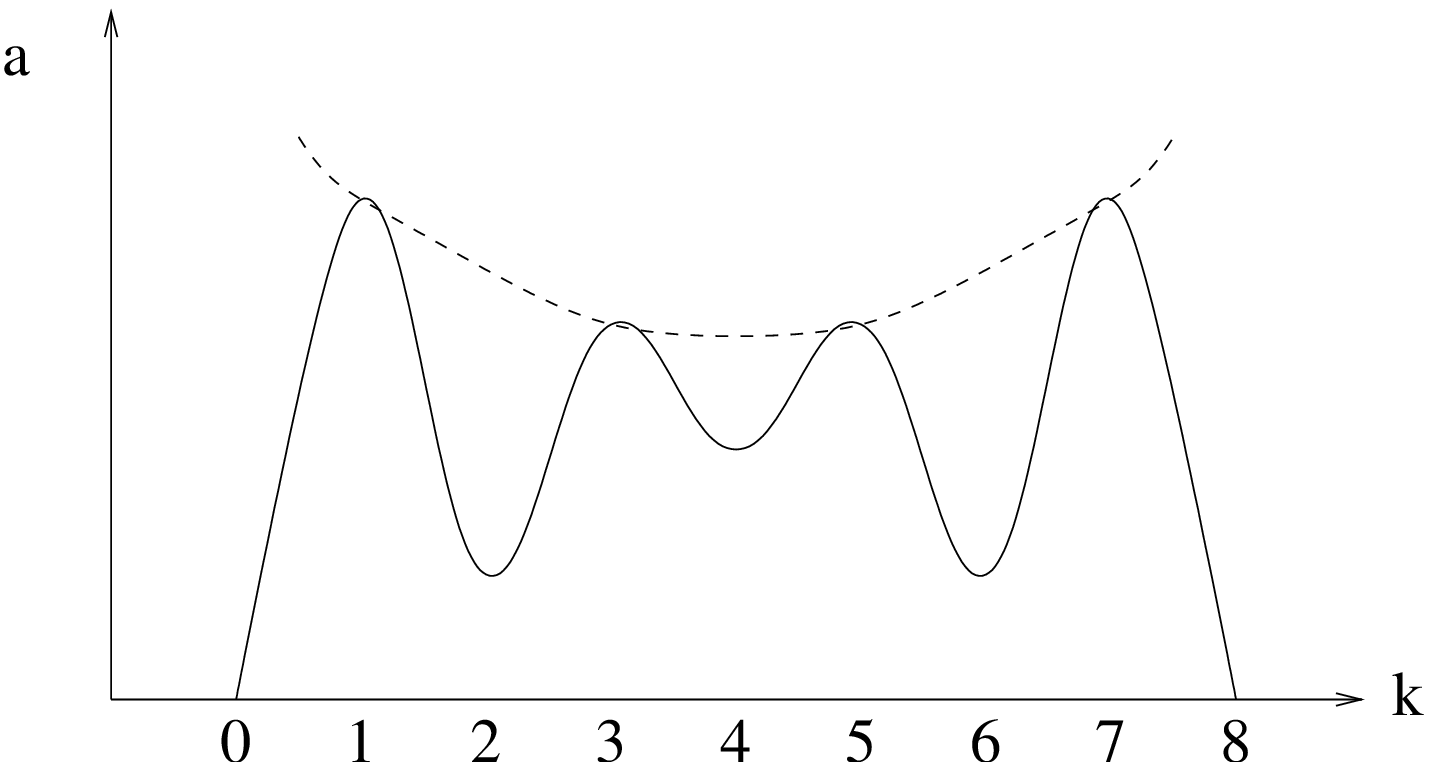}}
\end{center}
\vspace{.1in}
\end{minipage} 
\]
This was proved in~\cite[Theorem 5.4]{\emahajan}.

On the other hand,
if we think of the lists as $\cd$-monomials
then the pattern seems to become reverse unimodal.
This is the content of Conjecture~\ref{conjecture_unimodal}.
The analogy here is far from being clear.
We propose that the source of reverse unimodal behaviour
lies in the dotted line in the figure.

\section{A recursion for the {\bf cd}-index} \label{s:recur}

In this section,
we give a recursion for computing
the $\cd$-index of an Eulerian poset
in terms of certain polynomial sequences.

\subsection{Two polynomial sequences}

We use induction to define two homogeneous polynomial sequences
$\phi_m$ and $\phi_m^{\prime}$ for $m \geq 0$
in the variables $\cv$ and $\dv$.
Let
\[
\phi_0 = \cv, \quad
\phi_0^{\prime} = -2, \quad
\phi_{m+1} = \cv \phi_{m} + \dv \phi_{m}^{\prime}, \quad
\phi_{m+1}^{\prime} = (-2) \phi_{m} - \cv \phi_{m}^{\prime}.
\]
\vanish{
$\phi_0 = \cv$ and
$\phi_0^{\prime} = -2$ and
$\phi_{m+1} = \cv \phi_{m} + \dv \phi_{m}^{\prime}$
and
$\phi_{m+1}^{\prime} = (-2) \phi_{m} - \cv \phi_{m}^{\prime}$.
}
Note that the definitions are arranged so that
$(\av-\bv)^m \bv = \phi_m + \bv \phi_m^{\prime}$
holds for all $m \geq 0$.
This follows from the inductive definition and
the identities
$\av - \bv = \cv - 2 \bv$ and $(\av - \bv) \bv = \dv - \bv \cv$.

\subsection{An alternate definition}

An alternate definition of the 
{\bf ab}-index of a graded poset $P$
is given by assigning weights to each chain in $P$.
For a chain 
$c = \{\hz = x_0 < x_1 < \cdots < x_k = \ho \}$
define the {\it weight} of the chain to be the product
$\wt(c) = w_1 \cdots w_n$,
where
$$w_i = \left\{\begin{array}{l l}
             \bv, & \text{if} \ \ i \in \{\rho(x_1), \ldots, \rho(x_{k-1})\}, \\
             \av-\bv, & \text{otherwise}.
             \end{array}\right.$$

\noindent
Hence the weight of the chain is given by
\[
\wt(c) = (\av-\bv)^{\rho(x_0,x_1)-1} \bv (\av-\bv)^{\rho(x_1,x_2)-1} \bv
\cdots \bv (\av-\bv)^{\rho(x_{k-1},x_k)-1}.
\]
Then it follows from the definition that
the {\bf ab}-index of $P$ is given by
$\Psi(P) = \sum_c \wt(c)$,
where $c$ ranges over all chains in the poset $P$.
We rewrite this sum as follows.
\begin{equation} \label{eqn}
\Psi(P) = (\av-\bv)^n + \sum_{\hz < x < \ho} (\av-\bv)^{\rho(x)-1} \bv \Psi([x,\ho]).
\end{equation}
In other words,
we group together terms by the element of the smallest rank
in a chain.

We are primarily interested in the {\bf cd}-index.
So now we restrict ourselves to the class of Eulerian posets.
By definition,
every interval of an Eulerian poset is also an Eulerian poset.
Hence the terms 
$\Psi(P)$
and
$\Psi([x,\ho])$
are expressible in the variables $\cv$ and $\dv$.
So we can think of them as the {\bf cd}-index of the respective posets,
which agrees with our earlier notation.
Our goal is to write an expression for 
$\Psi(P)$
that involves only $\cv$ and $\dv$.

\subsection{The recursion}

Depending on whether the parity of $n$ is even or odd, 
we may write 
$(\av-\bv)^{n} = (\cv^2 - 2 \dv)^{n/2}$ or
$(\av-\bv)^{n} = \cv (\cv^2 - 2 \dv)^{n-1/2} - 2 \bv (\cv^2 - 2 \dv)^{n-1/2}$.
We first do the case when $P$ has odd rank.
Then we may write equation~\eqref{eqn} as
\[
\Psi(P) = (\av-\bv)^n + \sum_{\hz < x < \ho} (\phi_{\rho(x)-1} 
+ \bv \phi_{\rho(x)-1}^{\prime}) \Psi([x,\ho]) .
\]
Dropping all the terms that begin with a $\bv$, we obtain
\[
\Psi(P) = (\cv^2 - 2 \dv)^{\rho(P)-1/2} + 
\sum_{\hz < x < \ho} \phi_{\rho(x)-1} \Psi([x,\ho]).
\]
If $P$ has even rank 
then 
we replace the term $(\cv^2 - 2 \dv)^{\rho(P)-1/2}$ by
$\cv (\cv^2 - 2 \dv)^{\rho(P)-2/2}$.
This gives us a nice recursion for computing the {\bf cd}-index
of an Eulerian poset.
As a special case,
if $P$ is the Boolean lattice of odd rank
then we get
\[
\Psi(B_n) = (\cv^2 - 2 \dv)^{n-1/2} + 
\sum_{k=1}^{n-1} \binom {n}{k} \phi_{k-1} \Psi(B_{n-k}).
\]
We may write a similar statement
for $n$ even.

\section{The cubical lattice} \label{s:cube}

In this section,
we lay down the algebraic framework
to study the $\cd$-index of the cubical lattice.
The results will be cubical analogues 
of those obtained in Sections~\ref{s:algebra} and
\ref{s:duality}.
Let $C_{n+1}$ be the face lattice 
of the $n$ dimensional cube.

\subsection{The basic setup}
For $v$, a {\bf cd}-monomial of degree $n$,
let $\gamma(v)$ be the coefficient of 
$v$ in $\Psi(C_{n+1})$.
In more fancy language,
$\gamma(v) = \langle \delta_v,\Psi(C_{n+1}) \rangle$.
We then extend the definition to $\Ff$ by linearity.
Figure~\ref{f:U} shows the $\cd$-index of the cubical lattice
for small ranks.
An important distinction between the cubical and the
Boolean lattice is that  
$\gamma(v) \not = \gamma(v^*)$.

\begin{proposition}[Ehrenborg-Readdy] \label{l:H}
There is a well-defined linear map $H \colon \Ff \rightarrow \Ff$
given by the initial conditions
$$H(1) = 0, \: \: H(\cv) = 2 \dv, \: \: H(\dv) = \cv \dv + \dv \cv$$
and the rule
$H(u v) = H(u) v + u H(v)$, 
such that 
$$\Psi(C_{n+1}) = \Psi(C_n) \cv + H(\Psi(C_n)).$$
\end{proposition}

Let $\hat H \colon \Ff \rightarrow \Ff$
be the linear map defined by 
${\hat H}(u) = H(u) + u \cv$ for $u \in \Ff$.
The first few values are also shown in
Figure~\ref{f:val} in Section~\ref{s:algebra}.
Note that $\hat H$ is defined on $\Ff$ and not on $\Fhat$.
Also observe that the definition of $\hat H$ 
is arranged so that the equation
\begin{equation} \label{e:h}
{\hat H}(\Psi(C_n)) = \Psi(C_{n+1}) \quad \text{holds for} \quad n \geq 1.
\end{equation}
We want to view $\Ff$ as a comodule over $\Fhat$.
To that end,
define
$\delta : \Ff \rightarrow \Ff \tensor \Fhat$ by
$\delta(u) = \Delta(u) + u \tensor \ev$.

The analogue of Theorem~\ref{t:coderivation}
may now be stated as follows.

\begin{theorem} \label{t:cubical_coderivation}
Let $\delta$ and $\hat H$ be as defined above. Then
\[
\delta \circ {\hat H}(u) = 
(2 (\id \tensor {\hat G}) + {\hat H} \tensor \id) \circ \delta(u).
\] 
\end{theorem}

\begin{figure}
$$  
\begin{array}{r c l }
\Psi(C_1) & = & 1 \\
\Psi(C_2) & = & \cv \\
\Psi(C_3) & = & \cv^2 + 2 \dv \\
\Psi(C_4) & = & \cv^3 + 4 \cv \dv + 6 \dv \cv \\
\Psi(C_5) & = & \cv^4 + 6 \cv^2 \dv + 14 \dv \cv^2 +
                16 \cv \dv \cv + 20 \dv^2
\end{array} 
$$
\caption{The {\bf cd}-index of the Cubical lattice for 
          ranks 1,2,3,4 and 5.}
\label{f:U}
\end{figure}

This may be proved directly just as Theorem~\ref{t:coderivation}.
We will prove it by proving its dual version.
We mention that this equation looks unfamiliar 
and we have never encountered it before.

\subsection{The dual setup}

Dualise the maps 
$\hat H$ and $\delta$ to get
the corresponding dual maps
$\hat H^*$ and $\delta^*$. 
Using the identification of 
$\Ff^*$
with
$\Ff$,
we obtain a map, 
$T \colon \Ff \rightarrow \Ff$
and a module map
$\centerdot : \Ff \tensor \Fhat \rightarrow \Ff$.
We continue to denote the module map by $\centerdot$
because it is induced from the $\centerdot$ product
on $\Fhat$ via the inclusion map
$\Ff \into \Fhat$.
We may also note that the maps
$T$ and $\centerdot$ have degree $-1$.

Now we state the dual version of Theorem~\ref{t:cubical_coderivation}.
The proof will be given a little later.

\begin{theorem} \label{t:cubical_derivation}
We have
$T \circ \centerdot = \centerdot \circ (2 (\id \tensor S) + T \tensor \id).$ 
This may also be expressed as 
$T(u \centerdot v) = T(u) \centerdot v + 2 (u \centerdot S(v))$
for $u \in \Ff$ and $v \in \hat F$.
Also,
$T(1) = 0$.
\end{theorem}

We now state the dual version of the property
${\hat H}(\Psi(C_n)) = \Psi(C_{n+1})$
given by equation~\eqref{e:h}.

\begin{lemma} \label{l:cubical_reduction}
Let $v$ be any {\bf cd}-monomial of positive degree. Then
$\gamma(T(v)) = \gamma(v)$.
\end{lemma}

An explicit description of the map $T$ 
is straightforward to obtain and
is given as follows.

\begin{lemma} \label{l:definition_T}
Let $m_1,m_2,\ldots,m_k$
be non-negative integers. The map $T$ is given by
$
T(\cv^{m_1} \dv \cv^{m_2}  \dv \ldots \dv  \cv^{m_k}) =
\cv^{m_1-1} \ldots \dv  \cv^{m_i}  \dv \ldots \cv^{m_k} +\\
\sum_{i=2}^k 2 \cv^{m_1} \ldots \dv  \cv^{m_i-1}  \dv \ldots \cv^{m_k} +
\sum_{i=1}^{k-1} 2 \cv^{m_1} \ldots \dv  \cv^{m_i}  \cv  \cv^{m_{i+1}} \dv \ldots \cv^{m_k}.
$
\end{lemma}

\begin{lemma} \label{l:alternate_T}
Let $m_1,m_2,\ldots,m_k$
be non-negative integers. The map $T$ is given by
$
T(\cv^{m_1} \dv \cv^{m_2}  \dv \ldots \dv  \cv^{m_k}) =
2 S(\cv^{m_1} \dv \cv^{m_2}  \dv \ldots \dv  \cv^{m_k}) -
\cv^{m_1-1} \ldots \dv  \cv^{m_i}  \dv \ldots \cv^{m_k} 
$.
\end{lemma}

The above lemma follows by simply comparing 
the explicit descriptions of the maps $S$ and $T$
given by Lemmas~\ref{l:definition_S} and \ref{l:definition_T}
respectively.
This relation between $T$ and $S$
can be used to derive Theorem~\ref{t:cubical_derivation}
from Theorem~\ref{t:derivation}.

\medskip
\begin{proof_}{{\it Proof of Theorem~\ref{t:cubical_derivation}.}}
We use Theorem~\ref{t:derivation} and 
Lemma~\ref{l:alternate_T} and the theorem follows 
from the following sequence of equalities.
$$
\begin{array}{r c l}
T(u \centerdot v) & = & 2 S(u \centerdot v) - ((u \centerdot v) - 1) \\
                  & = & 2 S(u \centerdot v) - ((u - 1) \centerdot v) \\
                  & = & 2 S(u) \centerdot v + 2 u \centerdot S(v) - 
                                              ((u - 1) \centerdot v) \\
                  & = & ( 2 S(u) - (u - 1)) \centerdot v + 
                                              2 u \centerdot S(v) \\
                  & = &  T(u) \centerdot v + 2 (u \centerdot S(v)).
\end{array}
$$
If the first letter of $u$ is $\cv$ then
$u - 1$ refers to the $\cv\dv$-monomial obtained by deleting it
else it refers to the element 0.
\end{proof_}

\medskip
\subsection{Simple applications}

We conclude this section by writing analogues (without proof)
to the results of Sections~\ref{ss:sa} and \ref{ss:ssa}.
We will not consider analogues to the results
of Sections~\ref{s:ln}-\ref{s:balance} in this paper.
That would be a project in itself.

\begin{lemma} \label{l:ceasy}
Let $u$ and $v$ be $\cd$-monomials of degree $m$ and $n$ 
respectively. Then
\[
\gamma(u \centerdot v) = \binom {m+n+1}{m} 2^{n+1} \gamma(u) \beta(v).
\]
\end{lemma}

\begin{example}
Following the lines of the computation that we made for
$\beta(\underbrace {1 \centerdot 1 \centerdot \ldots \centerdot 1}_{n})$
in Example~\ref{eg:warmup} (Section~\ref{ss:sa}),
we see that 
$\gamma(\underbrace {1 \centerdot 1 \centerdot \ldots \centerdot 1}_{n}) = 
2^n n!.$
\end{example}

\begin{lemma} \label{l:cubical_add_on}
Let $u$ and $v$ be {\bf cd}-monomials of the same degree
and $w$ be any {\bf cd}-monomial. Then we have,
$$
\begin{array}{r c l}
\beta(u) > \beta(v) & \mbox{\rm iff} 
         & \gamma(w \centerdot u) > \gamma(w \centerdot v) \\
\gamma(u) > \gamma(v) & \mbox{\rm iff} 
         & \gamma(u \centerdot w) > \gamma(v \centerdot w).
\end{array}
$$
\end{lemma}

\begin{lemma} \label{l:cubical_switch}
Let $u,v$ and $w$ be {\bf cd}-monomials. Then we have,
$\gamma(u \centerdot v) = \gamma(u \centerdot v^*)$
and
$\gamma(u \centerdot v \centerdot w) = \gamma(u \centerdot w \centerdot v)$.
\end{lemma}

Recall that while proving this result for the Boolean lattice
the base case for induction was 
$\beta(\ev \centerdot v) = \beta(\ev \centerdot v^*)$.
For proving the above result, 
we use the base case
$\gamma(1 \centerdot v) = \gamma(1 \centerdot v^*)$.
Or we could also directly use Lemma~\ref{l:ceasy}.

\begin{lemma} \label{lemma_cubical_identities}
Let $u$, $v$ and $w$ be any {\bf cd}-monomials. Then we have,
$\gamma(u \dv \cv \dv v \dv) =
\gamma(u \dv \cv \dv v^* \dv)$
and 
$\gamma(u \dv \cv \dv 
v \dv \cv \dv w) =
\gamma(u \dv \cv \dv 
v^* \dv \cv \dv w).$
\end{lemma}

\begin{lemma} \label{lemma_cubical_cc_d}
Let $u$ and $v$ be {\bf cd}-monomials. Then
$\gamma(u \dv v) \geq \gamma(u \cv^2 v)$, 
with equality if $u$ and $v$ are both empty.
\end{lemma}

\section {More on the algebra $\Fhat$} 
\label{s:fun}

In Section~\ref{s:duality},
we defined an associative algebra structure on $\Fhat$.
In later sections, we used it effectively
to study the function $\beta$ that we were interested in.
In this section, we study this algebra in its own right.

\begin{theorem} \label{t:free}
Under the $\centerdot$ product, 
${\mathcal F}$ is a free algebra on countably
many generators. There are two natural sets of generators
$\{1,{\bf d},{\bf d}^{2},\cdots \}$ and 
$\{1,{\bf c}^{2},{\bf c}^{4},\cdots \}$.
\end{theorem}
\begin{proof}
We show the first part. 
The second part is left to the reader.
We prove the lemma in two steps.
In the first step, we show that 
$1,{\bf d},{\bf d}^{2},\cdots$ 
generate $\Ff$ 
and in the second step,
we show that they do not satisfy any relation.

\medskip
\noindent
{\em Step 1}:
We do a forward induction on the degree of the {\bf cd}-monomial
and for each degree, we do a backward induction 
on the number of $\dv$'s that it ends with.

Let $v$
be any {\bf cd}-monomial.
Write $v = u \cv^m \dv^k$, where $m > 0$
and $u$ ends in $\dv$.
By Lemma~\ref{l:bullet}, we obtain
$2 v = (u \cv^{m-1}) \centerdot \dv^k 
             - u \cv^{m-2} \dv^{k+1}$.
The monomial
$u \cv^{m-1}$
has a lower degree, while
$u \cv^{m-2} \dv^{k+1}$ 
has the same degree but 
ends with a larger number of $\dv$'s.
Therefore by our induction hypothesis,
these monomials
can be expressed in terms of our generators
and hence so can $v$.
This completes the induction step.

As an example, for $v = \cv^3$, write
$2 \cv^3 = \cv^2 \centerdot 1 - \cv \dv$.
Repeating the process on $\cv^2$ and $\cv \dv$,
we get
$4 \cv^3 = (\cv \centerdot 1 - \dv) \centerdot 1 - 1 \centerdot \dv$.
Substituting, $2 \cv = 1 \centerdot 1$, we get
$8 \cv^3 = 1 \centerdot 1 \centerdot 1 - 2 \dv \centerdot 1 
                 - 2 (1 \centerdot \dv)$.

\medskip
\noindent
{\em Step 2}:
Suppose there is a homogeneous relation between our generators, say
\[
1 \centerdot v_1 + \dv \centerdot v_2 + \cdots 
               + \dv^l \centerdot v_{n+1} = 0, \quad v_i \in \Ff.
\]
We first show that $v_1 = 0$.
Let $v_1 = \sum c_i w_i$, 
where $w_i$ are {\bf cd}-monomials of the same degree as $v_1$
and $c_i$ are constants.
By Lemma~\ref{l:bullet},
observe that the term $1 \centerdot w_i$, 
which occurs in $1 \centerdot v_1$, 
is (in general) a sum of two terms,
of which exactly one begins with a $\cv$.
This term does not appear in any other product term.
So, we conclude that $c_i = 0$, which says that $v_1 = 0$.
This reduces our relation to 
$\dv \centerdot v_2 + \cdots 
               + \dv^l \centerdot v_{n+1} = 0$, with $v_i \in \Ff$.
Repeating essentially the same argument, we get 
$v_i = 0$ for all $i$.
\end{proof}
\noindent
For an equivalent result, see~\cite[Theorem 3.4]{\bn}.

Next we recall the Dehn-Sommerville relations for
the flag $f$-vector of an Eulerian poset 
(Theorem~\ref{t:Dehn_Sommerville}) and
show that they are equivalent
to certain simple identities that exist in $\Fhat$.

\begin{theorem} \label{t:Dehn_Sommerville}
For an Eulerian poset $P$ of rank $n+1$ and 
a subset $S \subseteq [n]$,
if $\{i,k\} \subseteq S \cup \{0,n+1\}$ such that $i < k$,
and $S$ contains no $j$ such that $i < j < k$, then
\begin{equation} \label{e:Dehn_Sommerville}
\sum_{j=i+1}^{k-1} (-1)^{j-i-1} f_{S \cup j}^{n+1} (P) =
f_S^{n+1}(P)(1 - (-1)^{k-i}).
\end{equation}
\end{theorem}
\noindent
We begin with equation~\eqref{e:coalg_morph}
from Section~\ref{s:algebra}.
$$
\Delta(\Psi(P)) = \sum_{\hz < x < \ho} \Psi([\hz,x]) \tensor \Psi([x,\ho]).
$$
The dual to this equation is the identity
$$
\Psi_P^*(u \centerdot v) = \sum_{\hz < x < \ho} 
\Psi_{([\hz,x])}^*(u)  \Psi_{([x,\ho])}^*(v),
$$
where $\Psi_{P}^*(w)$
denotes the coefficient of $w$ in $\Psi(P)$.
Using the associativity of the $\centerdot$ product, 
we may write
$$
\Psi_P^*(u_1 \centerdot u_2 \ldots \centerdot u_{k+1}) = 
\sum_{\hz < x_1 < \cdots x_k < \ho} 
\Psi_{([\hz,x_1])}^*(u_1) \ldots \Psi_{([x_k,\ho])}^*(u_{k+1}).
$$
Recall that if $P$ is an Eulerian poset of rank $n+1$,
then $\Psi_{P}^*(\cv^n) = 1$.
Hence setting $u_i = \cv^{a_i}$, we get
$\Psi_P^*(\cv^{a_1} \centerdot \cv^{a_2} \ldots \centerdot \cv^{a_{k+1}}) = 
\sum 1$, 
where summation ranges over the set
$\{\hz < x_1 < \cdots x_k < \ho : \: \:
\rho(x_1) = a_1+1, \rho(x_2) = a_1+a_2+2, \ldots, 
\rho(x_k) = a_1+ \cdots + a_k + k\}$.
Here $\rho$ denotes the rank function of the poset.
From this observation, we conclude the following.

\begin{lemma} \label{l:relation}
We have
$
\Psi_P^*(\cv^{a_1} \centerdot \cv^{a_2} \centerdot 
\ldots \centerdot \cv^{a_{k+1}}) = f_S,
$
the component of the flag $f$-vector of $P$ for
$S = \{a_1+1, a_1+a_2+2,\ldots,a_1+ \cdots + a_k + k\}$.
\end{lemma}

\begin{lemma} \label{l:Euler}
Let $n$ be a positive integer. Then
\begin{equation} \label{e:sidentity}
(\cv^0) \centerdot (\cv^{n-1}) - (\cv^1) \centerdot (\cv^{n-2}) 
 + \cdots + (-1)^{n-1}(\cv^{n-1}) \centerdot (\cv^0) =
(1 + (-1)^{n+1}) (\cv^{n}).
\end{equation}
\end{lemma}

This follows directly from the definition of the
$\centerdot$ product
given by Lemma~\ref{l:bullet_list}.
Now let $P$ be any Eulerian poset of rank $n+1$. 
Applying $\Psi_P^*$ to both sides
of equation~\eqref{e:sidentity}
and applying Lemma~\ref{l:relation}, 
we obtain the Euler relation
\[
f_1 - f_2 + \cdots + (-1)^{n-1} f_{n} = (1 + (-1)^{n+1}).
\]
This corresponds to the case when
$S = \phi, i=0$ and $k = n+1$
in equation~\eqref{e:Dehn_Sommerville}.
To get the general case,
first rewrite the identity in 
Lemma~\ref{l:Euler} with
$n+1 = k-i$ to obtain
$$
\sum_{j=i+1}^{k-1} \cv^{j-i-1} \centerdot \cv^{k-j-1} = 
(1 - (-1)^{k-i}) \: \: \cv^{k-i-1}.
$$
Let $S = \{s_1,s_2,\ldots,s_l\}$.
For simplicity, we only explain the case when
$\{i,k\} \subseteq S$.
Let
$s_a = i$ and $s_{a+1} = k$.
Now pre and post multiply the above identity
by 
$\Pi_{i=1}^a \cv^{s_i - s_{i-1} -1}$ and
$\Pi_{i=a+2}^l \cv^{s_i - s_{i-1} -1}$ respectively.
Here the product $\Pi$ is taken
with respect to the $\centerdot$ product.

Now apply $\Psi_P^*$ to both sides
and use Lemma~\ref{l:relation}, 
to get equation~\eqref{e:Dehn_Sommerville}.

\subsection*{Acknowledgement}
This work was done in Spring 1998 under the supervision of Ken Brown. 
I would like to thank him for his insightful comments,
perusal of earlier versions of this paper and countless discussions. 
Midway through this project,
Harold Fox provided me numerical data till rank 14
for the Boolean and cubical lattice.
It verified the already established results and
pointed me to new ones.
He deserves a lot of credit for making this paper more complete.
I would also like to thank Richard Ehrenborg, 
who introduced me to this area of mathematics.
The motivation for this paper came from an earlier joint work with him.
I also thank Marcelo Aguiar and Sam Hsiao for many discussions.

%\bibliographystyle{hamsplain}
%\bibliography{boolean}

\begin{thebibliography}{10}

\bibitem{aguiar02:_infin_hopf}
M.~Aguiar, \emph{Infinitesimal {H}opf algebras and the cd-index of polytopes},
  Discrete Comput. Geom. \textbf{27} (2002), no.~1, 3--28, Geometric
  combinatorics (San Francisco, CA/Davis, CA, 2000).

\bibitem{aguiar00:_infin_hopf}
Marcelo Aguiar, \emph{Infinitesimal {H}opf algebras}, New trends in Hopf
  algebra theory (La Falda, 1999), Amer. Math. Soc., Providence, RI, 2000,
  pp.~1--29.

\bibitem{bayer85:_gener_dehn_sommer_euler}
Margaret~M. Bayer and Louis~J. Billera, \emph{Generalized {D}ehn-{S}ommerville
  relations for polytopes, spheres and {E}ulerian partially ordered sets},
  Invent. Math. \textbf{79} (1985), no.~1, 143--157.

\bibitem{bayer91}
Margaret~M. Bayer and Andrew Klapper, \emph{A new index for polytopes},
  Discrete Comput. Geom. \textbf{6} (1991), no.~1, 33--47.

\bibitem{billera00:_monot}
Louis~J. Billera and Richard Ehrenborg, \emph{Monotonicity of the cd-index for
  polytopes}, Math. Z. \textbf{233} (2000), no.~3, 421--441.

\bibitem{billera97}
Louis~J. Billera, Richard Ehrenborg, and Margaret Readdy, \emph{The
  $c$-$2d$-index of oriented matroids}, J. Combin. Theory Ser. A \textbf{80}
  (1997), no.~1, 79--105.

\bibitem{billera00:_noncom}
Louis~J. Billera and Niandong Liu, \emph{Noncommutative enumeration in graded
  posets}, J. Algebraic Combin. \textbf{12} (2000), no.~1, 7--24.

\bibitem{ehrenborg01:_euler}
Richard Ehrenborg, \emph{$k$-{E}ulerian posets}, Order \textbf{18} (2001),
  no.~3, 227--236.

\bibitem{ehrenborg02:_inequal}
\bysame, \emph{Inequalities for polytopes and zonotopes}, preprint, 2002.

\bibitem{ehrenborg:_inequal}
Richard Ehrenborg and Harold Fox, \emph{Inequalities for $\cd$-indices of joins
  and products of polytopes}, Combinatorica, to appear.

\bibitem{ehrenborg98:_maxim}
Richard Ehrenborg and Swapneel Mahajan, \emph{Maximizing the descent
  statistic}, Ann. Comb. \textbf{2} (1998), no.~2, 111--129.

\bibitem{ehrenborg96}
Richard Ehrenborg and Margaret Readdy, \emph{The $\mathbf r$-cubical lattice
  and a generalization of the ${\bf cd}$-index}, European J. Combin.
  \textbf{17} (1996), no.~8, 709--725.

\bibitem{ehrenborg98:_coprod}
\bysame, \emph{Coproducts and the $cd$-index}, J. Algebraic Combin. \textbf{8}
  (1998), no.~3, 273--299.

\bibitem{joni79:_coalg}
S.~A. Joni and G.-C. Rota, \emph{Coalgebras and bialgebras in combinatorics},
  Stud. Appl. Math. \textbf{61} (1979), no.~2, 93--139.

\bibitem{kalai88}
Gil Kalai, \emph{A new basis of polytopes}, J. Combin. Theory Ser. A
  \textbf{49} (1988), no.~2, 191--209.

\bibitem{purtill93:_andre}
Mark Purtill, \emph{Andr\'e permutations, lexicographic shellability and the
  $cd$-index of a convex polytope}, Trans. Amer. Math. Soc. \textbf{338}
  (1993), no.~1, 77--104.

\bibitem{reading02:_bases_euler}
Nathan Reading, \emph{Bases for the flag $f$-vectors of eulerian posets},
  preprint, 2002.

\bibitem{stanley89:_log}
Richard~P. Stanley, \emph{Log-concave and unimodal sequences in algebra,
  combinatorics, and geometry}, Graph theory and its applications: East and
  West (Jinan, 1986), New York Acad. Sci., New York, 1989, pp.~500--535.

\bibitem{stanley94:_flag}
\bysame, \emph{Flag $f$-vectors and the $cd$-index}, Math. Z. \textbf{216}
  (1994), no.~3, 483--499.

\bibitem{stenson02:_relat}
Cathy Stenson, \emph{Relationships among flag $f$-vector inequalities for
  polytopes}, preprint, 2002.

\bibitem{sundaram94:_cohen_macaul}
Sheila Sundaram, \emph{The homology representations of the symmetric group on
  {C}ohen-{M}acaulay subposets of the partition lattice}, Adv. Math.
  \textbf{104} (1994), no.~2, 225--296.

\bibitem{sundaram95}
\bysame, \emph{The homology of partitions with an even number of blocks}, J.
  Algebraic Combin. \textbf{4} (1995), no.~1, 69--92.

\end{thebibliography}
%input a.bbl
\providecommand{\bysame}{\leavevmode\hbox to3em{\hrulefill}\thinspace}

\end{document}